\documentclass{amsart}
\usepackage[british]{babel}
\usepackage{amsmath,amssymb,amsthm,amscd}
\usepackage{graphicx,subfigure}
\usepackage{color}
\usepackage{pb-diagram} 
\usepackage{tabularx}
\usepackage{multirow}
\usepackage{longtable}
\usepackage{hyperref}
\usepackage[foot]{amsaddr}

\usepackage{caption,rotating}
\usepackage{comment}

\newtheorem{theorem}{Theorem}[section]

\theoremstyle{remark} 
\newtheorem{remark}[theorem]{Remark}
\newtheorem{algorithm}[theorem]{Algorithm}

\numberwithin{equation}{section}

\newcommand{\te}{\theta}
\newcommand{\vp}{\varphi}

\newcommand{\cD}{{\mathcal D}}

\newcommand{\RR}{{\mathbb R}}
\newcommand{\CC}{{\mathbb C}}

\newcommand{\TT}{{\mathbb T}}
\newcommand{\ZZ}{{\mathbb Z}}
\newcommand{\NN}{{\mathbb N}}

\newcommand{\abs}[1]{|{#1}|}

\newcommand{\norm}[1]{\|{#1}\|}

\newcommand{\aver}[1]{\langle{#1}\rangle}
\newcommand{\Lie}[1]{\mathfrak{L}_{#1}}

\newcommand{\R}[1]{\mathfrak{R}_{#1}}
\newcommand{\Dioph}[2]{\cD_{#1,#2}}

\def\to{\rightarrow}

\def\epsilon{\varepsilon}

\def\ii{\mathrm{i}}
\def\ee{\mathrm{e}}

\def\dif{ {\mbox{\rm d}} }
\def\Dif{ {\mbox{\rm D}} }

\def\pd{ \partial }

\def\mani{{\mathcal{M}}}

\def\torus{{\mathcal{K}}}

\def\H{h}                   

\def\DeltaK{\Delta K}

\def\cteOmega{c_{\Omega,0}}
\def\cteDOmega{c_{\Omega,1}}

\def\cteJ{c_{J,0}}
\def\cteDJ{c_{J,1}}
\def\cteDDJ{c_{J,2}}

\def\cteJT{c_{J^{\top}\hspace{-0.5mm},0}}
\def\cteDJT{c_{J^{\top}\hspace{-0.5mm},1}}
\def\cteG{c_{G,0}}
\def\cteDG{c_{G,1}}
\def\cteDDG{c_{G,2}}

\allowdisplaybreaks

\begin{document}

\title[KAM theory for periodically or quasi-periodically systems]{
Constructive QP-Time-dependent KAM Algorithm for Lagrangian Tori}

\date{\today}

\author{Renato Calleja$^{\mbox{\textdagger} }$}
\address[\textdagger]{IIMAS, Universidad Nacional Auti\'onoma de M\'exico, Apdo. Postal 20-126, C.P. 01000, M\'exico D.F., M\'exico.}
\email{calleja@mym.iimas.unam.mx}

\author{Alex Haro$^{\mbox{\textdaggerdbl}, \mbox{\textasteriskcentered}}$}
\address[\textdaggerdbl]{Departament de Matem\`atiques i Inform\`atica, Universitat de Barcelona,
Gran Via 585, 08007 Barcelona, Spain.}
\address[\textasteriskcentered]{ Centre de Recerca Matem\`atica, Edifici C, Campus Bellaterra, 08193 Bellaterra, Spain.}
\email{alex@maia.ub.es}

\author{Pedro Porras$^{\mbox{\textdagger}, \mbox{\textdaggerdbl}}$}
\email{pedro.porras@iimas.unam.mx}
\thanks{Corresponding author: Pedro Porras (\texttt{pedro.porras@imass.unam.mx}).}

\begin{abstract}
In this paper, 
we present an algorithm to compute a fiberwise Lagrangian torus in quasi-periodic (QP) Hamiltonian systems, 
whose convergence is proved in the \cite{calleja2025constructive}. 
We exhibit the algorithm with two models.
The first is a Tokamak model \cite{chandre2005control,valvo2021hamiltonian}, 
which proposes a control method to create barriers to the diffusion of magnetic field lines through a small modification in the magnetic perturbation. 
The second model \cite{del2000self}, 
known as the vorticity defect model, 
describes the nonlinear evolution of localized vorticity perturbations in a constant vorticity flow. 
This model was originally derived in the context of plasma physics and fluid dynamics.
\end{abstract}

\maketitle

\tableofcontents

\section{Introduction}

One of the classical approaches for studying non-integrable Hamiltonian systems is through small perturbations of an integrable system expressed in action-angle coordinates. 
This approach led to the development of KAM theory in the 1950s and 1960s by Kolmogorov, 
Arnold, 
and Moser in \cite{kolmogorov1954conservation, Arnold63a, moser1962invariant}, 
(see \cite{celletti2006kam, dumas2014kam} for a brief history), 
aimed at demonstrating the persistence of invariant tori when perturbing integrable systems. 
However, 
the perturbative nature of KAM theory can pose challenges for practical applications, 
as perturbations are not always arbitrarily small. 
A relatively recent approach is the so-called parameterization method, 
which was first applied in the context of KAM theory in  \cite{de2005kam,de2001tutorial}. 
This method does not rely on the classical action-angle variables are not in the hypotheses of the theorem \cite{de2005kam} but in the consequences, 
and the method in the proof leads to an efficient algorithms of computation and computer-assisted proofs 
(CAPs) as in \cite{HaroCFLM16,FiguerasHL17}. 
Similarly,
\cite{valvo2021hamiltonian} discusses CAPs using normal forms.

These mathematical developments have found a natural application in the study of celestial mechanics, 
fluids and plasma transport models just to mention a few. 
In the context of plasma physics, 
non-autonomous Hamiltonian systems arise in Tokamak models,
devices that confine plasma in a toroidal geometry using magnetic fields. 
Developed in the 1950s by Tamm and Sakharov, 
Tokamaks represent one of the most promising approaches to achieving controlled nuclear fusion. 
The dynamics of the confined plasma can be considered as non-autonomous Hamiltonian systems, 
where quasi-periodic perturbations arise from variations in external fields. 
These systems challenge traditional perturbative techniques, highlighting the importance of modern tools like the parameterization method. 
For instance, the persistence of invariant tori in such models it is essential for maintaining stable plasma confinement, 
as demonstrated in \cite{mackay1984transport}. By bridging experimental plasma physics and advanced mathematical frameworks, 
recent research has provided valuable insights into the design and optimization of next-generation fusion reactors.

The stability of magnetic structures, 
such as quasi-periodic flux surfaces, 
is essential for maintaining confinement and preventing energy losses or plasma instabilities. 
KAM methods are employed to study the robustness of these surfaces, 
providing valuable insights into designing more efficient and stable confinement configurations \cite{boozer2004physics}. 
These developments demonstrate the broader impact of quasi-periodic Hamiltonians, 
as they not only enhance our understanding of intrinsic dynamics but also contribute to technological advancements in areas as diverse as space exploration and sustainable energy production.

Motivated by these needs, in \cite{calleja2025constructive}, 
we prove a KAM theorem,
using the parameterization method to look for an invariant torus in a QP-time-dependent Hamiltonian systems that depends periodically or quasi-periodically on time, 
with $\ell$ external frequencies.
The proof follows an a-posteriori approach, 
that uses a quasi-Newton method. 
This iterative method begins with an initial parameterization of an approximately invariant torus, that can be obtained in several ways, such as through
numerical calculations, Lindstedt series, the integrable case, etc,
meaning it approximately satisfies the invariance equation as defined by the invariance error. 
The approximation is refined
by applying corrections that reduce quadratically the invariance error.

In this paper we implement the method suggested in the proof of \cite{calleja2025constructive}, 
to compute invariant tori in diferent models,
which stand out for their relevance in plasma physics and fluid dynamics, both with structures governed by Hamiltonian systems and chaotic transport. 

The work \cite{chandre2005control}, 
fosuses on diffusion control in Hamiltonian systems with one and a half degrees of freedom, 
introducing a method to create barriers to the diffusion of magnetic field lines in a Tokamak time-periodic model. 
These barriers, 
represented by invariant tori, 
are generated through localized chaos control using a correction term of a smaller order than the original perturbation. 
The technique finds critical applications in magnetic fusion scenarios, 
such as restoring confinement in stochastic regions.

On the other hand, 
in \cite{del2000self}
a transported field and a velocity field are dynamically coupled via a mean field term coming from the advection-diffusion equation. 
Specifically, the analysis focuses on the single-wave model, 
which reduces to a quasi-periodic pendulum system characterized by the interaction of two external frequencies and one internal frequency, 
forming a three-dimensional torus.
This model involves self-consistent chaotic transport, 
offering a robust framework for studying the dynamics of coupled systems and their invariant structures.

The application of the KAM method for time-dependent periodic or quasi-periodic invariant tori is well-founded in the context of both scenarios, 
given their Hamiltonian structure and intrinsic time dependence. 
In the Tokamak case, 
KAM theory offers a rigorous framework to establish the persistence of invariant tori.
In the second case, 
the quasi-periodic nature emerging from the single-wave model aligns naturally with the assumptions of KAM theory. 
This method predicts the persistence of invariant tori while preserving the structured transport properties of the system. 
By leveraging the KAM framework, 
both models benefit from its robust capacity to explain the stability and persistence of invariant structures in time-dependent dynamical systems, 
providing deep insights into their behavior and stability.

The structure of this paper is as follows:
In Section ~\ref{sec:method}, 
we present the method of the KAM theorem, 
detailing the necessary geometric, 
analytical, 
and dynamical framework. 
This includes the formulation of quasiperiodic Hamiltonian systems, 
the invariance equation, 
and solutions to cohomological equations. 
Additionally, 
we discuss the KAM method as a quasi-Newton scheme, 
using the approximately invariant torus as the initial guess and refining successive approximations through corrections.
In Section ~\ref{sec:num:thm:Tokamak}, 
we apply the method to a Tokamak model with one and a half degrees of freedom. 
Using the Poincare section and Birkhoff averages, 
we compute internal frequencies,
using the tecniques in \cite{das2017quantitative},
and construct a two-dimensional approximately invariant torus as the initial guess. 
We then refine this torus using the KAM method and analyze its continuation with respect to parameters. 
Furthermore, 
we examine the behavior of Sobolev norms as indicators of potential breakdowns of the invariant torus, see \cite{calleja2010numerically, calleja2010computation}. 
In Section ~\ref{sec:num:thm:Pendulum}, 
we study one more system: 
a quasiperiodic pendulum with two external perturbations, with three-dimensional invariant tori. 
We present three scenarios;
quasi-periodically perturbed rotor (related to single wave model), 
quasi-periodically perturbed pendullum (rotational tori) and 
quasi-periodically perturbed pendullum (librational tori).
We emphasize the challenges of implementing these three-dimensional objects numerically, 
paying special attention to optimizing computational resources, 
such as RAM memory, and ensuring code efficiency.
Lastly, 
in Appendix~\ref{sec:appen:Tokamak}, 
we present a table for the Tokamak model with the parameter value $0.0004$. 
The table includes the time and memory usage percentages at each step of the quasi-Newton method iterations.
Through these sections, 
we demonstrate not only the versatility of the KAM method in various scenarios but also the importance of a robust numerical approach for tackling problems in high-dimensional quasiperiodic Hamiltonian systems. 
This work lays the groundwork for future research on the optimization and application of KAM methods in physical and mathematical contexts.
For instance, 
to compute invariant tori in the Restricted Three-Body Problem (RTBP), 
around the equilibrium points $L_4$ and $L_5$.
\section{KAM algorithms for Langrangian invariant tori}\label{sec:method}

In this section, 
we present the method of the KAM theorem, 
as described in \cite{calleja2025constructive}, 
introducing the necessary geometric, 
analytic, 
and dynamical framework. 
This includes the formulation of the time-quasiperiodic Hamiltonian, 
the associated Hamiltonian structure, 
the corresponding symplectic vector field, 
the invariance equation, 
and the formal solutions to the cohomological equations, 
considering Diophantine frequency vectors. 
With this foundation, 
we establish the KAM method and algorithm, 
which requires an approximately invariant torus that satisfies the invariance equation approximately, 
generating what we call the invariance error. 
This error serves as the initial guess for a quasi-Newton method, 
where the initial approximation is improved by adding a correction determined by a first-order approximation of the new invariance error, 
hence the name of the process. 
Subsequently, 
leveraging the symplectic structure, 
we apply a symplectic change of coordinates that reduces the linearized equations to cohomological equations. 
Finally, 
we describe methods and algorithms for the continuation of tori with respect to parameters, 
using the first variation with respect to these parameters as a key tool.

\subsection{Invariant Lagrangian tori in quasi-periodic Hamiltonian systems}\label{ssec:KAM:Method}
We denote by $\TT^n=\RR^n/\ZZ^n$ the $n$-dimensional torus.
Now, 
let us consider space is an open subset $\mani \subset \TT^{n_1}\times\RR^{n_2}\times\RR^{n}$, 
where $n_1 + n_2 =n$,
equipped with an exact symplectic form $\boldsymbol{\omega}$, 
a Riemannian metric $\boldsymbol{g}$, 
and a linear isomorphism $\boldsymbol{J}$,
that satisfy the relation $\boldsymbol{\omega}_z(\boldsymbol{J}_z u, v) = \boldsymbol{g}_z(u, v)$ for all $u, v \in T_z\mani$, 
and $z = (x,y)$. 
The corresponding matrix representations of these objects are given by $\Omega, G, J: \mani \to \RR^{2n \times 2n}$. 
More specifically, they fulfill the conditions:
\[
    \Omega^\top = -\Omega\,, \qquad G^\top = G\,, \qquad J = -\Omega^{-1} G\,.
\]
In this geometrical setting is referred to as \textbf{Cases II} in \cite{HaroCFLM16}, 
and the particular case in which, moreover,
$J^2 = -I_{2n}$ is referred to as \textbf{Case III}.

In this setting, 
a quasi-periodic Hamiltonian sysmtem is given by a function $\H: \mani \times \TT^\ell \to \RR$,
and the external frecuency vector $\alpha \in \RR^\ell$, for which the equations are,
\begin{equation}\label{eq:equations}
    \left\{
        \begin{split} \dot z &= Z_\H(z, \varphi)= \Omega(z)^{-1} (\Dif_{z} \H(z,\varphi))^\top \,, \\
            \dot \varphi &= \alpha \,.
     \end{split} 
     \right.
\end{equation}
The quasi-periodic Hamiltonian flow, 
\begin{equation}\label{eq:flow}
\begin{aligned}
\tilde{\Phi}_t : \mani \times \TT^{\ell} \times \RR & \longrightarrow \mani \times \TT^{\ell} \\
(z, \varphi, t) & \longmapsto  \tilde{\Phi}_{t}(z, \varphi) := \begin{pmatrix} \phi_{t}(z,\varphi) \\ \varphi + t \alpha \end{pmatrix}  \,,
\end{aligned}
\end{equation}
where the evolution operator $\phi_t(z,\varphi)$ satisfies
\[
\begin{split}
    \frac{\partial \phi_t}{\partial t}(z,\varphi) &= Z_\H(\phi_t(z,\varphi),\varphi + t \alpha)\,,\\
    \phi_0(z,\varphi) &=z \,.
\end{split}
\]

The goal is, to find a parameterization given a internal frecuency vector $\omega \in \RR^n$,
\[
\begin{array}{rcl}
    \tilde{K} : \TT^d \times \TT^{\ell} & \longrightarrow & \mani \times \TT^{\ell} \,, \\
(\theta, \varphi) & \longmapsto & \tilde{K}(\theta, \varphi) := \begin{pmatrix} K(\theta,\varphi) \\ \varphi \end{pmatrix}  \,,
\end{array}
\]
that satisfies the following invariance equation,
\begin{equation}
    \phi_{t}\left( K(\theta, \varphi), \varphi \right) = K(\theta + \omega t, \varphi + \alpha t)\,.
\label{eq:inv:flow}
\end{equation} 
or the infinitesimal version,
\begin{equation}
Z_{h}(K(\theta, \varphi), \varphi) = -\Lie{\omega, \alpha} K(\theta, \varphi)\,,
\label{eqn:invariance_equation}
\end{equation} 
where, $\Lie{\omega,\alpha}$ is defined in equation . 
The torus $\tilde \torus = \left\{(K(\theta,\varphi),\varphi) | (\theta,\varphi)\in \TT^n\times\TT^\ell \right\}$,
parameterized by $\tilde K$, is an invariant torus with frecuency vector $(\omega,\alpha)$ for the system \eqref{eq:equations}.
We assume that the parametrization is regular, meaning that $\Dif_\theta K(\theta,\varphi)$ has rank $n$.
The way the torus $\tilde \torus$ is embeded in $\mani \times \TT^\ell$ is given by the \textit{degree map} $\mathfrak{D} \in \ZZ^{n\times n_1}$ so that $K$ is of the form
\[
    K(\theta, \varphi) = \begin{pmatrix} \mathfrak{D} \theta \\ 0_{n_2} \\ 0_{n} \end{pmatrix} + K_p(\theta, \varphi) = \begin{pmatrix} \mathfrak{D} \theta \\ 0_{n_2} \\ 0_n\end{pmatrix} + \begin{pmatrix} K_x^p(\theta, \varphi) \\  K_y^p(\theta, \varphi)  \end{pmatrix} \,,
\]
where $K_p : \TT^n\times\TT^\ell \to \RR^{2n}$ is 1-periodic in $(\theta,\varphi)$, meaning $K_p(\theta + e_1,\varphi + e_2) = K_p(\theta,\varphi)$ for all $(e_1,e_2) \in \ZZ^n\times\ZZ^\ell$.

We introduce a tangent frame 
$L:\TT^n \times \TT^{\ell} \rightarrow \RR^{2n \times n}$,
and the normal frame $N:\TT^n \times \TT^{\ell} \rightarrow \RR^{2n \times n}$ 
is constructed as follows:
\begin{equation}\label{eq:L}
    L(\theta,\varphi):= \Dif_\theta K(\theta,\varphi)\,,
\end{equation}
and
\begin{equation}\label{eq:N}
N(\theta,\varphi):= L(\theta,\varphi) A(\theta,\varphi) + \tilde N(\theta,\varphi)\,,
\end{equation}
where,
\begin{align}
\tilde N(\theta,\varphi)={} & J(K(\theta,\varphi)) L(\theta,\varphi)B(\theta,\varphi)\,, \label{eq:tildeN} \\
B(\theta,\varphi)={}& (L(\theta,\varphi)^\top G(K(\theta,\varphi)) L(\theta,\varphi))^{-1}\,, \label{eq:B} \\
A(\theta, \varphi)={}& 
\begin{cases} -\displaystyle\frac{1}{2} (\tilde N(\theta,\varphi)^\top \Omega(K(\theta,\varphi))  \tilde N(\theta,\varphi)), & 
\text{if \textbf{Case II;}} \\
O_n\,, & \text{if \textbf{Case III.}} 
\end{cases}
\, \label{eq:A} 
\end{align}
The torsion matrix $T:\TT^n \times \TT^{\ell} \to\RR^{n\times n}$, 
given by 
\begin{align}
T(\theta, \varphi)={}& 
\begin{cases} - \dfrac{1}{2}\tilde N(\te,\vp)^\top \left(T_h(K(\te,\vp),\vp) +  T_h(K(\te,\vp),\vp)^\top\right)\tilde N(\te,\vp) & \text{if \textbf{Case II;}} \\ 
  \quad + \tilde N(\te,\vp)^\top T_h(K(\te,\vp),\vp)^\top N(\te,\vp)  & \\ 
  \quad + N(\te,\vp)^\top T_h(K(\te,\vp),\vp)\tilde N(\te,\vp),  & \\ 
  & \\ 
N(\te,\vp)^\top T_h(K(\te,\vp),\vp) N(\te,\vp)\,, & \text{if \textbf{Case III.}} 
\end{cases}
\, \label{eq:T} 
\end{align}
Moreover, 
we have that, 
given a parameterization of a torus 
$\torus=K(\TT^n \times \TT^{\ell})$
with, $T_\H:\mani \times \TT^{\ell} \to\RR^{2n\times 2n}$
\begin{align}
T_\H(z, \varphi)={}& 
\begin{cases} \Omega(z) \Big( \Dif_z Z_\H(z, \vp) & \text{if \textbf{Case II;}} \\ 
  \quad -  \Dif_z J(z) \left[Z_\H(z, \vp) \right] J(z)^{-1}& \\ 
    \quad - J(z)\Dif_z Z_\H(z, \vp) J(z)^{-1} \Big), & \\ 
  & \\ 
\Omega(z) \Big( \Dif_z Z_\H(z, \vp)  & \text{if \textbf{Case III.}} \\
\quad + \Dif_z J(z) \left[Z_\H(z, \vp) \right] J(z)& \\
    \quad + J(z) \Dif_z Z_\H(z, \vp) J(z) \Big)\,. &
\end{cases}
\, \label{eq:Th} 
\end{align}
\begin{remark}
    Additionally, 
    the symplectic form is canonical, 
    therefore,
    $\Omega=\Omega_0$, $G = I_{2n}$, and $J=\Omega_0$.
    Hence, 
    we have $\cteOmega=1$, 
    $\cteDOmega=0$,  
    $\cteG=1$, 
    $\cteDG=0$,  
    $\cteDDG=0$, 
    $\cteJ=1$, 
    $\cteDJ=0$, 
    $\cteDDJ=0$, 
    $\cteJT=1$, 
    and $\cteDJT=0$. 
    so the torsion is 
    \begin{equation}\label{eq:torsion:standar}
        T(\theta,\varphi) = N(\theta,\varphi)^\top T_\H(K(\theta, \varphi), \varphi) N(\theta,\varphi) \,.
    \end{equation}
    with
    \begin{equation}\label{eq:torsionH:standar}
        T_\H(K(\theta,\varphi),\varphi) = \Omega_0 \Dif_z Z_\H(K(\theta,\varphi), \vp) - \Dif_z Z_\H(K(\theta, \varphi), \vp) \Omega_0 \,.
    \end{equation}
In terms of the Hamiltonian, 
the torsion is expressed as follows
    \begin{equation}\label{eq:torsion:hamil}
            T(\theta,\varphi) = B(\theta,\varphi)^{-1} L(\theta,\varphi)^\top 
            S_\H(K(\theta,\varphi),\varphi) L(\theta,\varphi)B(\theta,\varphi)^{-1} \,.
    \end{equation}
with
\begin{equation}\label{eq:torsionS:hamil}
            S_\H(K(\theta,\varphi),\varphi) = \begin{pmatrix}
            \left( \Dif_{yy} - \Dif_{xx}\right)\H(K(\theta,\varphi),\varphi) & -2\Dif_{xy} \H(K(\theta,\varphi),\varphi) \\
            -2\Dif_{xy} \H(K(\theta,\varphi),\varphi)&  \left( \Dif_{xx} - \Dif_{yy}\right)\H(K(\theta,\varphi),\varphi) 
        \end{pmatrix} \,.
\end{equation}

\end{remark}
If the torus is approximately invariant, 
we define the \emph{error of invariance}
$E:\TT^n \times \TT^{\ell} \rightarrow \mani$ given by
\begin{equation}\label{eq:inv:err}
E(\theta, \varphi) = Z_{h}(K(\theta, \varphi), \varphi) + \Lie{\omega, \alpha} (K(\theta, \varphi)).
\end{equation}
The set of Diophantine vectors is defined as
\begin{equation}\label{eq:def:Dioph}
\Dioph{\gamma}{\tau} =
\left\{
(\omega, \alpha) \in \RR^d\times\RR^{\ell} \, : \,
\abs{k_1 \cdot \omega + k_2 \cdot \alpha} \geq \frac{\gamma}{|(k_1,k_2)|_1^{\tau}}
\,, 
\forall (k_1,k_2) \in\ZZ^d\times\ZZ^{\ell}\backslash\{(0,0)\}
\right\}\,.
\end{equation}
An essential condition in this theorem is the assumption that the frequency vector 
$(\omega,\alpha)$ for specific $\gamma >0$ and $\tau \geq n + \ell -1$, 
satisfies the Diophantine conditions.
The quantity $|(k_1,k_2)|_1$ represents the sum of the absolute values of each component of the vector $(k_1,k_2)$.
\begin{remark}
    It is not possible to represent exact Diophantine vectors numerically, 
    because of the limited precision of floating-point arithmetic. 
    This limitation means that we cannot accurately define irrational numbers, 
    and therefore true Diophantine vectors. 
    Instead, 
    we can approximate these vectors using the best possible numbers within the arithmetic precision we are working with. 
    For each component of the frequency vector that is not an integer, 
    we express it as a continued fraction. 
    We replace the final part of the continued fraction with ones to improve the approximation and make the components more accurate.
    In this way, 
    we could compute numerically Diophantine constants as in \cite{FiguerasHL17}.
\end{remark}
Finally, 
let us present notation regarding the cohomological equations which are central to KAM theory. 
This notation describes the relationship between a periodic function, 
$v:\TT^n\times\TT^{\ell} \rightarrow \RR$, 
and the frequency vector, 
$(\omega,\alpha) \in \RR^{d}\times \RR^{\ell}$. For this, 
we consider,
the Fourier expansion of a periodic function as
\[
    v(\theta,\varphi)=\sum_{k_1 \in \ZZ^d } \sum_{k_2 \in \ZZ^{\ell} } \hat v_{k_1,k_2} \ee^{2\pi \ii (k_1 \cdot \theta + k_2 \cdot \varphi)},
\]
\[
    \hat v_{k_1,k_2} =
    \int_{\TT^{\ell}} \int_{\TT^d} v(\theta, \varphi) \ee^{-2 \pi \ii (k_1 \cdot \theta + k_2 \cdot \varphi)} \dif \theta \, \dif \varphi \,,
\]
and introduce the notation $\aver{v}:=\hat v_{0,0}$ for the average. 
Hence,
the cohomological equation is the following,
\begin{equation}\label{eq:calL}
\Lie{\omega,\alpha} u = v- \aver{v}\,, \qquad
\Lie{\omega,\alpha} := -\left( \sum_{i=1}^d \omega_i \frac{\pd}{\pd \theta_i} +\sum_{j=1}^{\ell} \alpha_j \frac{\pd}{\pd \varphi_j} \right).
\end{equation}
The formal solution of equation \eqref{eq:calL}, 
with zero average, 
can be expressed as 
\begin{equation}\label{eq:small:formal}
\R{\omega,\alpha} v(\theta,\varphi) =\displaystyle \sum_{(k_1,k_2) \in \ZZ^n \times \ZZ^\ell \backslash \{(0,0)\} } \hat u_{k_1,k_2} \ee^{2\pi
\ii (k_1 \cdot  \theta + k_2 \cdot \varphi)}, \quad \hat u_{k_1,k_2} = \frac{-\hat
v_{k_1,k_2}}{2\pi \ii  (k_1 \cdot \omega + k_2 \cdot \alpha )}\,,
\end{equation}
with $(k_1, k_2) \neq (0,0)$ and $\hat u_{0,0}$ is a degree of freedom.
For $r \in \RR^+$, 
the Sobolev space $H^r$ is a Banach space consisting of functions $v$. 
These functions are characterized by the norm  

\begin{equation}
    \|v\|_r = \sqrt{\displaystyle \sum_{k_1 \in \ZZ^n} \displaystyle\sum_{k _2\in \ZZ^{\ell}} |2\pi\ii (k_1 + k_2)|^{2r} |\hat{v}_{k_1,k_2}|^2} < \infty.
\end{equation}

Here, 
$| \cdot |$ denotes the maximum norm in the spaces $\RR^n$ and $\CC^n$, 
meaning that for a vector $x = (x_1, \dots, x_n) \in \CC^n$, 
it is defined as  

\begin{equation}
|x| := \max_{j=1, \dots,n} |x_j|\,,
\end{equation}
similarly, 
this notation extends to real or complex matrices of arbitrary dimension, 
referring to the matrix norm induced by the corresponding vector norm.\\

The behavior of this norm is important in the study of the breakdown of quasi-periodic solutions. 
As discussed in \cite{calleja2010numerically, calleja2010computation},
the Sobolev norm serves as a key diagnostic tool for identifying the loss of smoothness in invariant tori. 
A rapid growth in this norm indicates the onset of breakdown, 
linking it to the regularity properties of the system. 
These ideas provide a numerical and theoretical foundation for detecting critical thresholds in quasi-periodic structures.

\subsection{Numerical method for torus computation}\label{sec:num:methd}
The constructuction of quasi-Newton method in \cite{calleja2025constructive} suggests to a numerical method. 
The approach involves starting with an approximately invariant torus, $K(\theta, \varphi)$,
meaning one with a non-zero invariance error $E$. The method is as follows. We obtain a new parameterization of the torus as 
        \[\bar K(\theta,\varphi) = K(\theta,\varphi) + \Delta K(\theta,\varphi)\,.\]
    The new invariance error is given by 
     \[
        \bar E(\theta, \varphi) = Z_H(\bar K(\theta, \varphi), \varphi) +\mathfrak{L}_{\omega, \alpha} \bar{K}(\theta,\varphi)\,.\\
    \]
        Expanding last equation, we have,
        \[
            \bar E(\theta, \varphi)= {\rm D}_z Z_H(K(\theta, \varphi), \varphi) \Delta K(\theta, \varphi) +\mathfrak{L}_{\omega, \alpha} \bar{K}(\theta,\varphi) + E(\theta, \varphi)
                +\Delta^2_{\theta} Z(\theta, \varphi)\,.
        \]
    To determine $\DeltaK$,
        we retain only terms up to first order, and solve the following equation,
        \[{\rm D}_zZ_H\bigr( K(\theta, \varphi), \varphi \bigr) \Delta K(\theta,\varphi) +\mathfrak{L}_{\omega, \alpha} \Delta K(\theta,\varphi) = -E(\theta, \varphi)\,.\] 
    To tackle this equation, we use the approximately symplectic framework represented by 
        \begin{equation}\label{eq:P}
            P(\theta, \varphi)= \left( L(\theta, \varphi) \, N(\theta, \varphi) \right) \,,
        \end{equation}
        such that, for some $\xi:\TT^n \times \TT^\ell \to \RR^{2n}$
        \[\Delta K(\theta, \varphi)= P(\theta, \varphi) \xi(\theta, \varphi)\,, \]
        i.e., $\xi(\theta, \varphi)$ is the new unknown.\\
    Leveraging certain geometric properties, we obtain,
        \begin{equation*} \label{eq:lin:1Enew}
        \begin{split}
            \left(\Lambda(\theta,\varphi)+ E_{red}(\theta,\varphi) \right) \xi(\theta,\varphi)
            + {} & \left( I_{2n} - \Omega_0 E_{sym}(\theta,\varphi)  \right)  \mathfrak{L}_{\omega, \alpha} \xi(\theta,\varphi) \\
            = {} &
            \Omega_0 P(\theta,\varphi)^\top \Omega(K(\theta,\varphi)) E(\theta,\varphi)
            \,,
        \end{split}
        \end{equation*}

    \begin{equation}\label{eq:Lambda}
        \Lambda(\theta,\varphi)
        = \begin{pmatrix}
            O_n &  T(\theta,\varphi) \\
            O_n  & O_n
        \end{pmatrix}\,.
    \end{equation}

    Neglecting quadratic terms, the solution to the above equation is approximated by solving a triangular system that requires handling cohomological equations,
         \begin{equation} \label{eq:lin:2Enew}
        \begin{split}
            \Lambda(\theta,\varphi) \xi(\theta,\varphi)
            + {} & \mathfrak{L}_{\omega, \alpha} \xi(\theta,\varphi) \\
            = {} &
            \Omega_0 P(\theta,\varphi)^\top \Omega(K(\theta,\varphi)) E(\theta,\varphi) \\
            = {} & \begin{pmatrix} - N(\theta,\varphi)^\top \Omega(K(\theta,\varphi)) E(\theta,\varphi) \\ 
            L(\theta,\varphi)^\top \Omega(K(\theta,\varphi)) E(\theta,\varphi)  
            \end{pmatrix} =:
    \begin{pmatrix}                                                                                                  
    \eta^L(\te,\vp) \\
    \eta^N(\te,\vp)
    \end{pmatrix}           \,.
        \end{split}
        \end{equation}
    To solve these equations, $\det \langle T\rangle ^{-1}$ must be different from zero, where,
             \begin{align}
                 \xi^N(\theta,\varphi)={} & \xi^N_{0,0} + \mathfrak{R}_{\omega,\alpha}(\eta^N(\theta,\varphi)) \,, \label{eq:etaN}\\
            \xi^L(\theta,\varphi)={} & \xi^L_{0,0} + \mathfrak{R}_{\omega,\alpha}(\eta^L(\theta,\varphi) - T(\theta,\varphi)
            \xi^N(\theta,\varphi))  \label{eq:etaL} \,,
        \end{align}
        with
        \[
        \xi^N_{0,0}= \langle T\rangle ^{-1} \langle \eta^L- T \mathfrak{R}_{\omega,\alpha}(\eta^N) \rangle \,.
        \]
If the torus is approximately invariant.
Therefore, the algorihtm is following:
\begin{algorithm}[Computation of the Adapted Frame and Torsion]\label{alg:P}
    Let $K$ satisfy the equation approximately. Compute the adapted frame $P$ and the reduced dynamics by following these steps:
\begin{enumerate}
    \item Derive $L(\theta, \varphi)$ using \eqref{eq:L}.
    \item Obtain $B(\theta, \varphi)$ from \eqref{eq:B}.
    \item Calculate $A(\theta, \varphi)$ based on \eqref{eq:A}.
    \item Evaluate $\tilde{N}(\theta, \varphi)$ according to \eqref{eq:tildeN}.
    \item Determine $N(\theta, \varphi)$ from \eqref{eq:N}.
    \item Compute $T_\H(\theta, \varphi)$ from \eqref{eq:Th}.
    \item Derive $T(\theta, \varphi)$ based on \eqref{eq:T}.
    \item Construct the frame $P(\theta, \varphi)$ using equation \eqref{eq:P}.
\end{enumerate}
\end{algorithm}

\begin{algorithm}[Correction of the Generating Torus]\label{alg:barK}
Let $K$ approximately satisfy the equation. Obtain the corrected generating torus by following these steps:

\begin{enumerate}
    \item Compute $P(\theta,\varphi)$  using Algorithm~\ref{alg:P}.
    \item Derive $E(\theta,\varphi)$ using \eqref{eq:inv:err}.
    \item Evaluate $\eta(\theta,\varphi)$, the right-hand side of \eqref{eq:lin:2Enew}.
    \item Solve the cohomological equation for $\xi^N$ and $\xi^L$, from \eqref{eq:etaN} and \eqref{eq:etaL} respectively, verifying $\aver{\eta^N}= 0$.
    \item Improve $K(\theta, \varphi)\leftarrow K(\theta,\varphi) + P(\theta,\varphi) \xi(\theta,\varphi)$.
\end{enumerate}
\end{algorithm}

\subsection{Numerical method for torus continuation}\label{sec:num:cont}
    Due to the nature of the method, 
    an initial torus is required. 
    Therefore, 
    if the system depends on a family of parameters, 
    we can use a previously computed torus for the next value of the parameter, 
    meaning that we can perform a continuation with respect to the parameter.
    In the following, 
    we describe a method for performing the continuation by considering the variational equation in the parameter.
    We know that the torus depends on the parameter, 
    i.e., 
    $K_{\epsilon}(\theta, \varphi)$,
    the first variational derivative in $\epsilon$ is
    \begin{equation}\label{eq:variational}
        K_{\epsilon + \delta_\epsilon }(\theta, \varphi) =  K_{\epsilon}(\theta, \varphi) + \frac{\partial}{\partial \epsilon} K_{\epsilon}(\theta, \varphi) \delta_\epsilon + \mathcal{O}(\delta_\epsilon^2)\,.
    \end{equation}
    Therefore, 
    we need to calculate $\frac{\partial}{\partial \epsilon} K_{\epsilon}(\theta, \varphi)$, 
    invoking the invariance equation \eqref{eqn:invariance_equation} and defined 
    $\Delta K_\epsilon(\theta,\varphi) := \frac{\partial}{\partial \epsilon}  K_\epsilon(\theta, \varphi) $, we obtain,
    \begin{equation}\label{eq:DeltaKeps}
        \Dif_z Z_\H ( K_{\epsilon}(\theta, \varphi),\varphi )\Delta K_{\epsilon}(\theta, \varphi) + \Lie{\omega,\alpha} \Delta K_\epsilon(\theta,\varphi) = - \frac{\partial}{\partial \epsilon} Z_\H ( K_{\epsilon}(\theta, \varphi), \varphi)\,,
    \end{equation}
    we proceed in a manner analogous to the \textit{numerical method}, 
    that is,
    we introduce a new $\xi_\epsilon:\TT^n \times \TT^{\ell} \to \RR^{2n}$, 
    such that $\Delta K_\epsilon(\theta,\varphi) = P(\theta, \varphi)\xi_\epsilon(\theta, \varphi)$,
    to obtain
\begin{equation} \label{eq:lin:KEps}
        \begin{split}
            \Lambda(\theta,\varphi) \xi_\epsilon(\theta,\varphi)
            + {} & \mathfrak{L}_{\omega, \alpha} \xi_\epsilon(\theta,\varphi) \\
            = {} &
            \Omega_0 P(\theta,\varphi)^\top \Omega(K_\epsilon(\theta,\varphi)) \frac{\partial}{\partial \epsilon} Z_{\H_\epsilon} ( K_{\epsilon}(\theta, \varphi), \varphi) \\
            = {} & \begin{pmatrix} - N(\theta,\varphi)^\top \Omega(K_\epsilon(\theta,\varphi)) \frac{\partial}{\partial \epsilon} Z_{\H_\epsilon} ( K_{\epsilon}(\theta, \varphi), \varphi) \\ 
                L(\theta,\varphi)^\top \Omega(K_\epsilon(\theta,\varphi)) \frac{\partial}{\partial \epsilon} Z_{\H_\epsilon} ( K_{\epsilon}(\theta, \varphi), \varphi)
            \end{pmatrix} =:
    \begin{pmatrix}                                                                                                  
    \eta_\epsilon^L(\te,\vp) \\
    \eta_\epsilon^N(\te,\vp)
    \end{pmatrix}           \,.
        \end{split}
        \end{equation}
    Therefore, 
    to find $\Delta K_\epsilon$, 
    it is enough to solve a cohomological equation and $\aver{\eta^N_\epsilon}=O_n$. 
    For the latter, 
    we have that
    \begin{equation} \label{eq:lin:KEps}
        \begin{split}
            \eta^N_\epsilon(\theta,\varphi) = {} &
            L(\theta,\varphi)^\top \Omega(K_\epsilon(\theta,\varphi)) \frac{\partial Z_{\H_\epsilon} }{\partial \epsilon}(K_\epsilon(\theta,\varphi),\varphi) \\
            = {} &
            \left(\Dif_\theta \left(\frac{\partial \H_\epsilon}{\partial \epsilon} (K_{\epsilon}(\theta, \varphi),\varphi)  \right) \right)^\top \,.
        \end{split}
    \end{equation}
    Hence,
    $\eta^N_\epsilon$ is the derivative of periodic functions, 
    and we obtain $\aver{\eta^N_\epsilon} = O_n$.
    \begin{algorithm}[Computation of the First Order Correction]\label{alg:DeltaKeps}
        Let $K_{\epsilon}(\theta,\varphi)$ satisfy the invariance equation approximately.
        To perform the continuation with respect to $\epsilon$, 
        compute the first variation $\Delta K_\epsilon$ by following these steps:
        \begin{enumerate}
            \item Calculate the derivative of $Z_{\H_\epsilon}(K(\theta, \varphi),\varphi)$ with respect to $\epsilon$.
            \item Evaluate $\eta_\epsilon(\theta, \varphi)$, the right-hand side of \eqref{eq:lin:KEps}.
            \item Solve the cohomological equation for $\xi^N_\epsilon$ and $\xi^L_\epsilon$.
            \item Compute $P(\theta,\varphi)$ using Algorithm~\ref{alg:P}.
            \item Improve $K_{\epsilon + \delta \epsilon}(\theta, \varphi)\leftarrow K_{\epsilon}(\theta,\varphi) + P(\theta,\varphi) \xi_\epsilon(\theta,\varphi)$.
            \item Save the $K_{\epsilon + \delta \epsilon}(\theta, \varphi)$ to be used in the continuation step.
        \end{enumerate}
    \end{algorithm}

    \begin{algorithm}[Continuation of the Torus family]\label{alg:ContinuationK}
        Let $K$ approximately satisfy the invariance equation. 
        Follow these steps to do a continuation respect to $\epsilon$:
    \begin{enumerate}
        \item Construct an approximate invariant parametrization.
        \item Improve the parametrization, $K_\epsilon(\theta, \varphi)$, using Algorithm~\ref{alg:barK}.
        \item Save the updated torus $K_\epsilon(\theta, \varphi)$.
        \item Define the continuation step $\delta$.
        \item Compute $K_{\epsilon + \delta \epsilon}(\theta, \varphi)$ using Algorithm~\ref{alg:DeltaKeps}.
        \item Repeat from step (2), until the desired value of $\epsilon$ is attained.
    \end{enumerate}
    \end{algorithm}
\begin{remark}
    The algorithms require to performe basic operations $(+,-,/,*)$, 
    compositions, 
    derivatives, and computing the $\R{\omega,\alpha}$ operator. 
    They do not require numerical integrations as in \cite{fernandez2024flow} and 
    can be applied as long as the family of tori exists given we have enough computational resources available.
\end{remark}
\begin{remark}
    Another approach to performing continuation with respect to parameters is the pseudo-arclength continuation method, 
    as described in \cite{gonzalez2022efficient}. 
    This work also provides a detailed implementation of the predictor used in the natural continuation method.
\end{remark}
\subsection{Comments on the implementations}\label{sec:comments}
In this section, 
we present the numerical implementation of the algorithms derived in the previous sections. 
While the theorem provides a rigorous mathematical foundation, 
its proof leads to an algorithm whose numerical implementation requires careful consideration of discretization
differentiation, equation-solving, 
and other numerical tools. 
To achieve this,
we have defined three key objects: 
\textbf{Matrix},
\textbf{Grid},
and \textbf{Fourier}, as in \cite{HaroCFLM16}
each playing a distinct role in the overall implementation.

Numerically, 
for the parameterization, 
we consider a set of sample points on a regular grid of size
$N_\theta=\left(N_{\theta,1},\dots, N_{\theta,n}\right)\in\NN^n$
and $N_\varphi=\left(N_{\varphi,1},\dots, N_{\varphi,\ell}\right)\in\NN^\ell$, 
such that,
\[
    \begin{split}
        \theta_{i}:= & \left(\theta_{i_1},\dots,\theta_{i_n}\right) = \left(\dfrac{i_1}{N_{\theta,1}},\dots,\dfrac{i_n}{N_{\theta,n}}\right)\,, \\
        \varphi_{j}:= & \left(\varphi_{j_1},\dots,\varphi_{j_\ell}\right) = \left(\dfrac{j_1}{N_{\varphi,1}},\dots,\dfrac{j_\ell}{N_{\varphi,\ell}}\right)\,,
    \end{split}
\]
where $i = (i_1,\dots,i_n)$, 
$j = (j_1,\dots,j_\ell)$, 
with $0\leq i_r \leq N_{\theta,r}\,$, 
$0\leq j_s \leq N_{\varphi,s}\,$,
further $1\leq r \leq n$ and $1\leq s \leq \ell$.
The total number of points is $N_D= N_{\theta, 1}\cdots N_{\theta, n}N_{\varphi, 1}\cdots N_{\varphi, \ell}$,
where $N_{\theta,i} = 2^{q_i}$, $N_{\varphi,j} = 2^{q_j}$,
with $q_i \in \NN$, $q_j\in\NN$ for $i=1,\dots,n$ and $j=1,\dots,\ell$.

We use the \textbf{Matrix} object, for functions, like; 
$K,E,\eta, \xi:\TT^{n}\times\TT^\ell\to\RR^{2n}$ or like; 
$L,N,T :\TT^{n}\times\TT^\ell\to\RR^{2n}$.
We can use the \textbf{Grid} or \textbf{Fourier} objects for function components like; 
$K_x, K_y :\TT^{n}\times\TT^\ell\to\RR^{n}$ or $\xi^N, \xi^L :\TT^{n}\times\TT^\ell\to\RR^{n}$. 
The \textbf{Grid} object is useful and efficients to do multiplications, sum, etc, or evaluating the vector field,
and we employ the \textbf{Fourier} object to calculate the derivative, 
or solving cohomological equations.
We switch between both representations with the aid of the Fast Fourier Transform (FFT),
which has a complexity of $\mathcal{O}(N_D \log(N_D))$.

One advantage of defining the torsion as in \eqref{eq:T},
see \cite{calleja2025constructive} for details, 
instead of $T(\theta,\varphi) = N(\theta,\varphi)^\top\Omega(K(\theta,\varphi))\left(\Dif_z Z(K(\theta,\varphi),\varphi) + \Lie{\omega,\alpha}K(\theta,\varphi)\right)$
is that we avoid the need to compute derivatives in Fourier space, as \cite{haro2019posteriori}. 
In other words, 
we do not generate a \textbf{Fourier} object, 
which would require additional memory and computational time for the derivative calculation. 
Instead, 
it suffices to perform grid-based multiplications.

\section{Numerical study of KAM Theory in a Tokamak model}\label{sec:num:thm:Tokamak}

In this section, 
we apply the KAM algorithm to a Tokamak model with one and a half degrees of freedom, 
following the background established in \cite{chandre2005control,valvo2021hamiltonian}. 
As in,
\cite{abdullaev2006construction}
we calculate the Poincaré section in the $(\psi,\theta)$ plane. 
From this section, and using Birkhoff averages, 
we obtain the internal frequency with the methods developed in \cite{das2017quantitative}.
With the internal frequency and the Poincare section, 
we impose the invariance equation to construct a two-dimensional approximately invariant torus with an external frequency equal to $1$. 
This torus serves as the initial guess for the method, 
to which we apply the KAM procedure.

Once the method converges, 
we perform a continuation with respect to the parameter and observe the behavior of some Sobolev norms. 
This provides estimates for the breakdown of analyticity of the invariant torus.

\subsection{A Tokamak model}\label{ssec:Tokamak:modl}
Tokamaks are devices for confining plasma in nuclear fusion,
see \cite{artsimovich1972tokamak}. 
They exhibit complicated dynamics that can be described using periodic non-autonomous Hamiltonian systems. 
In these systems, 
charged particles move along paths determined by magnetic fields, 
where the topological geometry is essential for maintaining plasma stability. 
KAM theory offers a theoretical framework for understanding how quasi-periodic orbits persist despite disturbances, 
thus contributing to the stability of invariant tori in phase space and controlling the plasma. 

Reference \cite{abdullaev2006construction} presents a non-autonomous Hamiltonian system that models the magnetic field lines of a Tokamak for plasma confinement, 
similar to equation \eqref{eq:Hamiltonians-Tokamak},
but without the term $\epsilon ^2$. 
Where, $\varphi$ represents the toroidal angle acting as the time, 
$\psi$ denotes the normalized toroidal flux, 
and $H$ stands for the poloidal flux. 
The poloidal angle $\theta$ is the conjugate variable to the action $H$.

In \cite{chandre2005control}, 
the authors apply, 
to the same model, 
a method they developed, 
introducing a control technique aimed at creating barriers between two chaotic regions by employing a control term proportional to $\epsilon^2$,
as depicted in figure (\ref{fig:strobos}).
We can observe that in \cite{valvo2021hamiltonian}, 
the authors use the same model and, 
through frequency analysis combined with a rigorous (Computer-Assisted), 
demonstrate that in the phase space of the magnetic field, 
the control term generates a set of invariant tori that act as transport barriers.
\begin{align}
    H(\theta, \psi; \varphi) &= \displaystyle \int \dfrac{ d\psi }{q(\psi)}  + \epsilon H_{1}(\theta, \varphi) + \epsilon^2 f_{2}(\theta, \varphi), \label{eq:Hamiltonians-Tokamak}\\
    q(\psi) &= \dfrac{4}{(2 -\psi)(2 - 2\psi + \psi^2)}, \nonumber \\
    H_{1}(\theta, \varphi) &= \cos(2\theta - \varphi) + \cos(3\theta - 2\varphi),\nonumber\\
    f_{2}(\theta,\varphi) &= \left( -\dfrac{1}{2} \dfrac{d}{d \psi} \left( \dfrac{1}{q(\psi)} \right) \right)\bigg|_{\psi_0} \left( \dfrac{2\cos(2\theta - \varphi)}{2w -1} + \dfrac{3\cos(3\theta - 2\varphi)}{3w -2} \right)^2,\nonumber \\
    w &= \dfrac{1}{q(\psi_0)}, \nonumber \\
    q_0 &= 0.35\,, \nonumber
\end{align}
where $\psi_0$ is the initial normalized toroidal flux on the Poincare section.
\begin{remark}
    In this case, 
    the geometric objects are given by
    $\Omega = \Omega_0, G = I_{2n}$ and $J=\Omega_0$.
    The dimension are $n_1=1, n_2=1, \ell=1$, and the degree map is $\mathfrak{D}=1$.
    Then, 
    the torus takes the form
\[
    \tilde K(\theta,\varphi) = \begin{pmatrix} K(\theta, \varphi) \\ \varphi \end{pmatrix} =
            \begin{pmatrix} \theta \\ 0 \\ 0 \end{pmatrix} +
            \begin{pmatrix} K_{x}^p(\theta, \varphi)  \\ K_{y}^p(\theta,\varphi) \\ \varphi \end{pmatrix}
                \,.
\]
\end{remark}
\subsection{Invariant torus}\label{ssec:inv:torus}
\begin{figure}[h!]
\centering         
\includegraphics[scale=0.8]{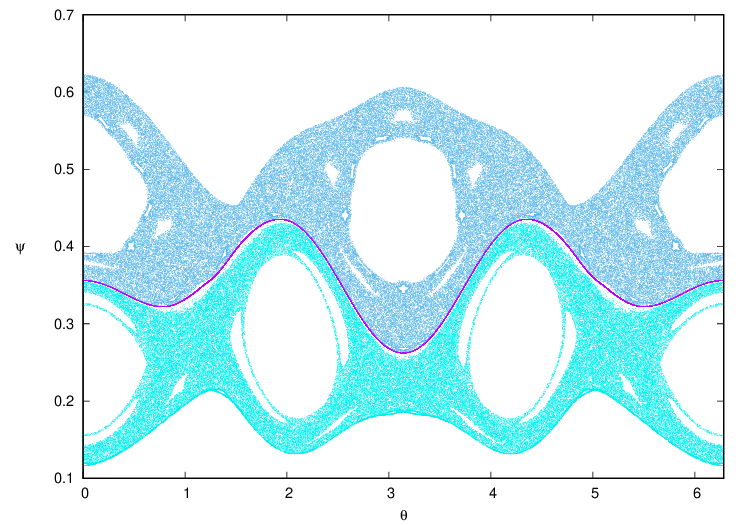}
\caption{The colors Cyan, Magenta and Teal  correspond to $\varphi = 0.35, 0.3566248878338341665, 0.36$, respectively.}
\label{fig:strobos}
\end{figure}

As previously mentioned in this work and in \cite{calleja2025constructive}, 
the method requires the construction of an initial torus 
with a small invariance error. 
To achieve this, 
we calculate numerically the stroboscopic map,
figure (\ref{fig:strobos}),
from the model described in equation \eqref{eq:Hamiltonians-Tokamak},
corresponding to $\epsilon = 0.004$.
To compute the internal frequency $\omega \approx 0.57981245427252670451$, 
we used Birkhoff averages on stroboscopic map, 
following \cite{das2017quantitative}. 
The external frequency $\alpha$, 
corresponding to the time-dependent part, 
is $1$. 
On the other hand,
we construct the two-dimensional initial torus from this curve, 
corresponding to $K(\theta, 0)$, 
denoted by $K_0(\theta)$. 
To do so, 
we construct the initial guess,
we define $(\bar \theta_i, \bar \varphi_j)$ such that,
$\theta_{i} = \bar \theta_{i} + \omega t$ and $\varphi_{j} =\bar \varphi_{j} + \alpha t$,
and we impose that $K(\bar \theta_{i}, \bar \varphi_{j})$ parametrized an invariant torus, 
i.e., the invariance equation is satisfied, 
equation \eqref{eq:inv:flow},
such that,
\[
\Phi_{t}\left( K(\theta_{i} -\omega t, \varphi_{j} - \alpha t), \varphi_{j} - \alpha t \right) = K\left(\theta_{i}, \varphi_{j}  \right)\,,
\]
since we want to construct the torus from $K_{0}(\theta)$, 
we impose $\varphi_{j} - \alpha t = 0$, 
with $t^*=\frac{\varphi_{j_1}}{\alpha_1}$ such that solve the last equation,
so we obtain,
\begin{equation}\label{eq:inv-eq-init-guess}
    \Phi_{t^*}\left( K_{0} \left( \theta_{i} - \omega t^* \right), 0\right) = K\left(\theta_{i}, \varphi_{j}  \right)\,.
\end{equation}
In other words, 
to obtain 
$K\left(\theta_{i}, \varphi_{j}  \right)$ it is enough to integrate the flow for a time 
$t^*$ with initial condition $K_0(\theta_{i} - \omega t^*)$. 
To evaluate the initial condition required to satisfy \eqref{eq:inv-eq-init-guess}, we evaluate $K_0$ using splines.\\

\begin{figure}[htbp]
    \centering
    \begin{minipage}{0.95\textwidth}                                                                                                                                         
        \centering
        \includegraphics[width=\textwidth]{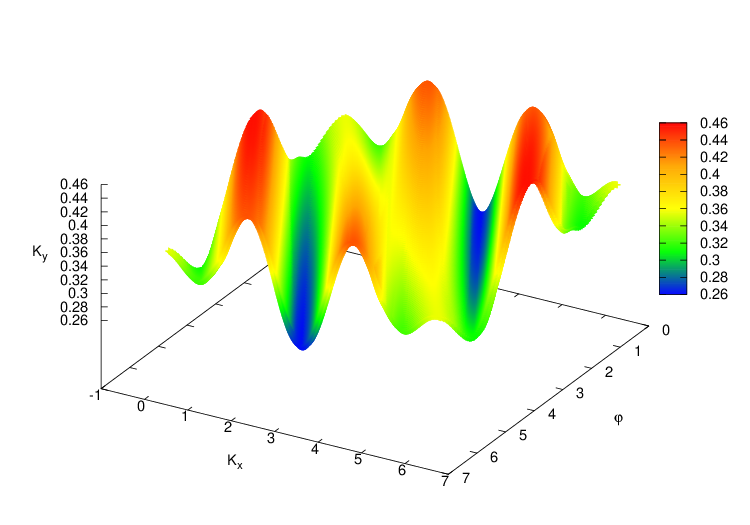}
        \caption{Parametrized surface of the invariant torus, 
        displaying the components $K_x(\theta,\varphi)$, 
        $K_y(\theta,\varphi)$, 
        and $\varphi$, 
        represented with $2^{23}$ Fourier modes and an invariance error equal to $3.96005\times10^{-16}$.
        Where blue represents the lowest height and red the highest, based on the $K_y(\theta,\varphi)$ component.
        }
    \label{fig:KAMtorusS}
    \end{minipage}
    \vspace{5mm} 
    \begin{minipage}{0.95\textwidth}
        \centering
        \includegraphics[width=\textwidth]{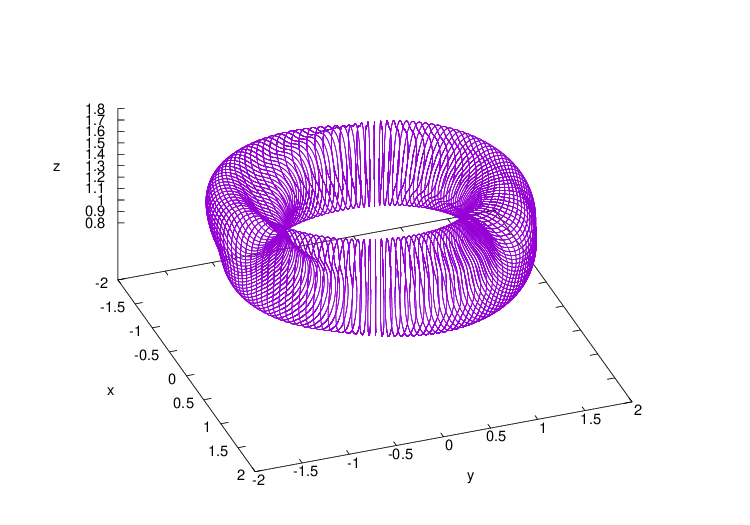}
    \end{minipage}
    \caption{Invariant torus in Cartesian coordinates.}
    \vspace{5mm}
    \label{fig:KAMtorus}
\end{figure}

We use a regular grid of points $N_\theta=(2^{n_\theta})$ and  $N_\varphi = (2^{n_\varphi})$,
with 
\[
        (\theta_i, \varphi_j ) = (\theta_{i_1}, \varphi_{j_1}) =  \left( \dfrac{i_1}{2^{n_\theta}}, \dfrac{j_1}{2^{n_\varphi}} \right)\,,
\]
with $1\leq i_1 \leq 2^{n_\theta}$ and $1\leq j_1 \leq 2^{n_\varphi}$.
Once we construct the invariant torus, 
with an initial invariance error approximate of $9.6866 \times 10^{-2}$ and $2^9 \times 2^9$, 
Fourier coefficients, 
it will serve as the initial guess of the quasi-Newton method.
We obtained an approximately invariant torus with an invariance error of $3.96005\times10^{-16}$, 
figures (\ref{fig:KAMtorusS}) and (\ref{fig:KAMtorus}),
and $2^{12}\times2^{11}$ Fourier coefficients in the $\theta$ and $\varphi$ direction, 
respectively.

\begin{figure}[h!]
\centering         
\includegraphics[scale=0.8]{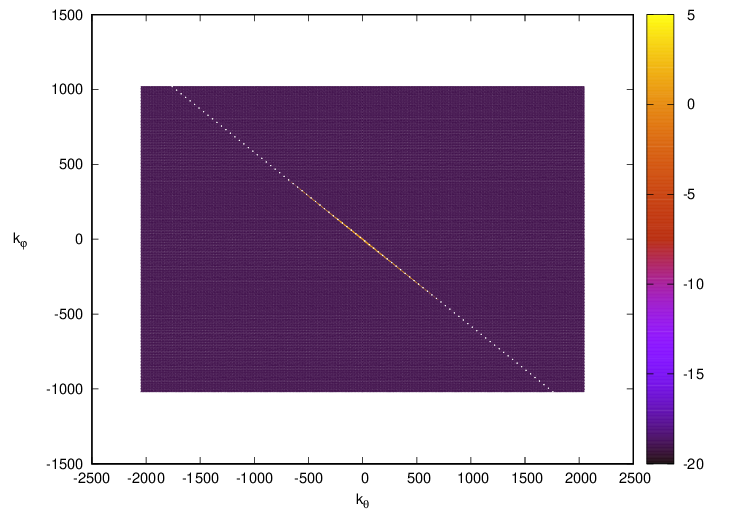}
\caption{Norm of the Fourier coefficients on a logarithmic scale. 
    The coefficients with the largest contributions, 
    the 'significant' ones, 
    align with the straight line shown in yelow in the image, 
    which fits the equation $y= -0.58013x$.}
\label{fig:norms:coefs}
\end{figure}
\begin{figure}[h!]
\centering         
\includegraphics[scale=0.8]{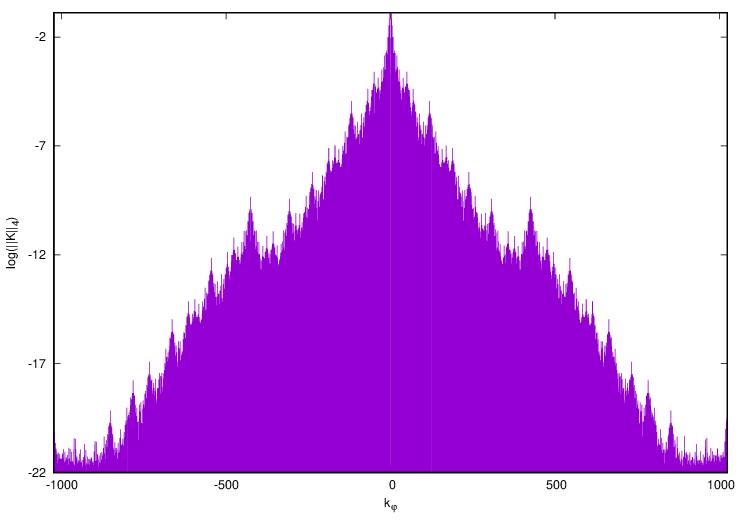}
\caption{For each $k_\varphi$, 
    we take the logarithm of the maximum of the $H^4$ Sobolev norm,
    equation \eqref{eq:H4-norm},
    of the Fourier coefficients for the torus with $\epsilon=0.004$.}
\label{fig:norms:coefs-max:400}
\end{figure}

\clearpage
\subsection{Continuation of the invariant torus}\label{ssec:inv:tori}
As we saw in Section ~\ref{ssec:inv:torus}, 
it is possible to apply the method without the need to perform a continuation over the perturbations starting from the integrable system. 
However, 
the strength of our approach lies in its ability to also perform a continuation over the perturbation parameter, 
either starting from the integrable system or from one already computed, as in Section~\ref{ssec:inv:torus}. 
To achieve this, 
we employ Algorithm~\ref{alg:ContinuationK}, 
using the case $\varepsilon = 0.004$ as the initial guess. 
The construction of the two-dimensional object follows an analogous process to that described in Section~\ref{sec:num:cont}. 
By performing the continuation as outlined in Algorithm~\ref{alg:DeltaKeps}, 
i.e., 
with a prediction for the following torus, 
we can compute the subsequent torus in four iterations of the Quasi-Newton step, 
compared to the six iterations required when improving the prediction.

\begin{figure}[h!]
\centering         
\includegraphics[scale=0.8]{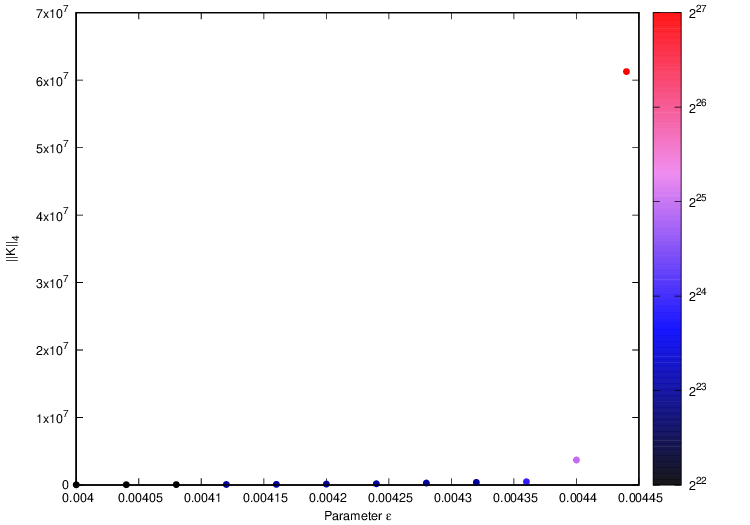}
    \caption{We observed that near the parameter value $\epsilon = 0.00445$, 
    there is an exponential growth, indicating signs of the breakdown of the torus.}
\label{fig:cont:coefs}
\end{figure}

\begin{figure}[h!]
\centering      
\includegraphics[scale=0.8]{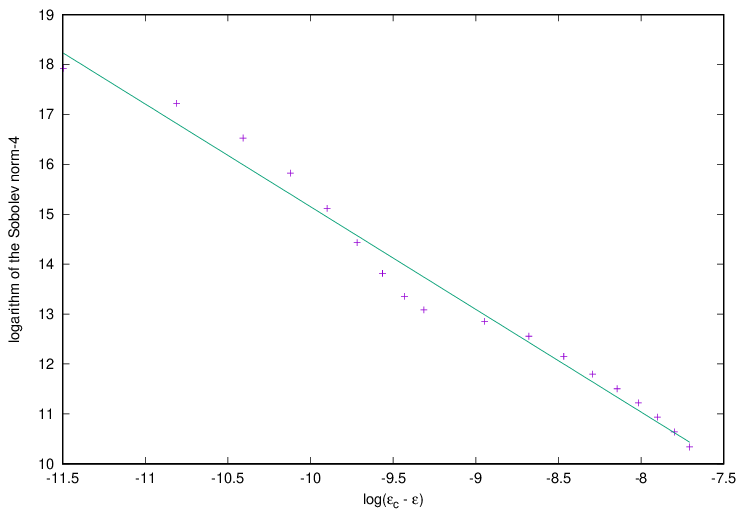}
    \caption{
			The logarithm of the $H^4$ Sobolev norm aligns with the linear fit $\log(\norm{H}_4) =-2.35233844056405\log(\epsilon_c -\epsilon) - 7.75624914299206$,
			with a correlation coefficient of $r=-0.987679863894756$ and a critical $\epsilon_c=0.0044615$.}
\label{fig:cont:sobonorms}
\end{figure}

While this paper aims to exemplify the suggested algorithms from \cite{calleja2025constructive}, 
we also seek to show the method's effectiveness by estimating the breakdown. 
We, now analyze the behavior of the of the $H^4$ Sobolev norm of $K$ defined as 
\begin{equation}\label{eq:H4-norm}
    \norm{K}_4 := \sqrt{\norm{K_x^p}_4^2 + \norm{K_y^p}_4^2 }\,.
\end{equation}
We observe, 
that in figure~\ref{fig:cont:coefs} Sobolev norm has a exponetial behavior \cite{calleja2010numerically}, 
we obtain an estimate for the blow-up exponent (see equation (12) in \cite{calleja2010computation}) in figure \ref{fig:cont:sobonorms}.
The fit to estimate the exponent also shows an oscilatory behavior in the Sobolev norm, figure \ref{fig:cont:sobonorms}.
These oscilations are not predicted by renormalization group, 
which is similar to the spin-orbit problem of Celestial Mechanics, 
\cite{calleja2024accurate}.

\begin{figure}[h!]
\centering         
\includegraphics[scale=0.8]{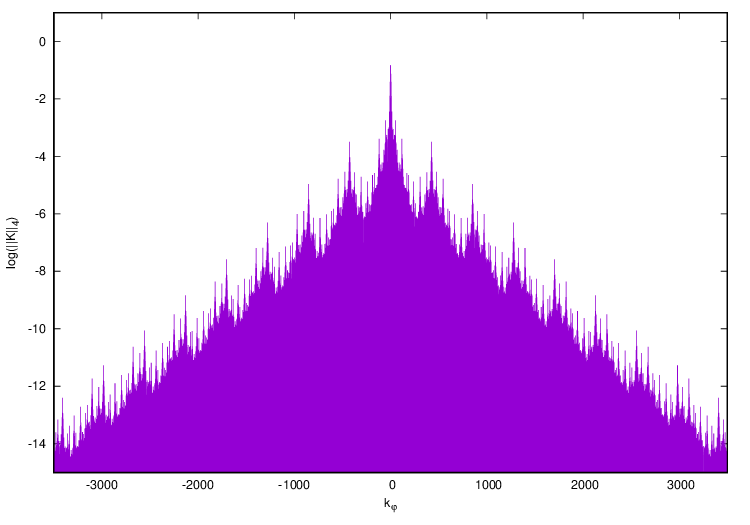}
\caption{For each $k_\varphi$, 
    we take the logarithm of the maximum of the H4-norm,
    equation \eqref{eq:H4-norm},
    of the Fourier coefficients for the torus with $\epsilon=0.00445$.}
\label{fig:norms:coefs-max:445}
\end{figure}
\clearpage
\section{Numerical study of KAM Theory in fluid and plasma transport models}\label{sec:num:thm:Pendulum}
In this section, we present the following model of a quasi-periodically forced pendulum with two driving frequencies,
\begin{equation}\label{eq:Transport:mod}
    \H=\dfrac{p^2}{2} +\epsilon_1\cos(q) - (\epsilon_2 + \epsilon_3\cos(\alpha_1 t)) \cos(q -\alpha_2t) \,,
\end{equation}
where, 
$\epsilon_1$,$\epsilon_2$ and $\epsilon_3$ are parameters,
with
\(\alpha_1 = \sqrt{3}\), and 
\(\alpha_2 = \frac{1+\sqrt{5}}{2}\). 
This pendulum model is motivated by the work in \cite{del2000self}. 
The main objective of this section is to demonstrate that, 
using the KAM algorithm, 
it is possible to perform a continuation with respect to parameters in the pendulum. 
To illustrate this, 
we present three distinct examples.
\begin{enumerate}
    \item 
    The \textit{the single-wave model} simplifies the vorticity mixing problem, 
    ultimately leading to the stream function. 
    A comprehensive derivation is available in References  \cite{del1998nonlinear,del1998weakly}, which provide a detailed analysis for interested readers.
    The continuation is carried out by initially constructing a planar pendulum as an initial guess and then gradually adapting it to the solitary wave model, 
    i.e.,
    $\epsilon_1=\epsilon_2=\epsilon_3=0$ to $\epsilon_1=0$, $\epsilon_2\neq0$ and $\epsilon_3\neq0$.
    \item Rotational tori. In this example, 
        we compute invariant tori corresponding to the libration tori of the pendulum, 
        starting from rotationonal torus of the simple pendulum, 
        and use continuation to transform it into the quasi-periodic case.
        We can construct an initial guess torus by setting $\epsilon_1\neq0$ and $\epsilon_2=\epsilon_3=0$, 
        and then continue varying the parameter $\epsilon_2$ and $\epsilon_3$.
    \item Lastly, 
        Librational tori we perform a similar analysis for the rotational tori of the simple pendulum. 
        This case is particularly interesting since the topology of librational tori differs from the others, 
        warranting detailed analysis.
        In this case, 
        the continuation is analogous to the rotational case.
\end{enumerate}
Unlike the previous section, 
these three-dimensional objects require significant effort in numerical implementation. 
Special attention must be paid to optimizing resources, 
such as RAM usage, 
as well as to code efficiency. 
These considerations are critical for ensuring accurate and efficient computations during the continuation of the tori. 
Through these examples, 
we not only demonstrate the capability of the KAM method to handle complex quasiperiodic systems but also highlight the importance of a robust numerical implementation.
\subsection{The single-wave model}\label{ssec:single-wave}
\begin{figure}[htbp]
    \begin{minipage}{0.45\textwidth}
        \centering
        \includegraphics[width=\textwidth]{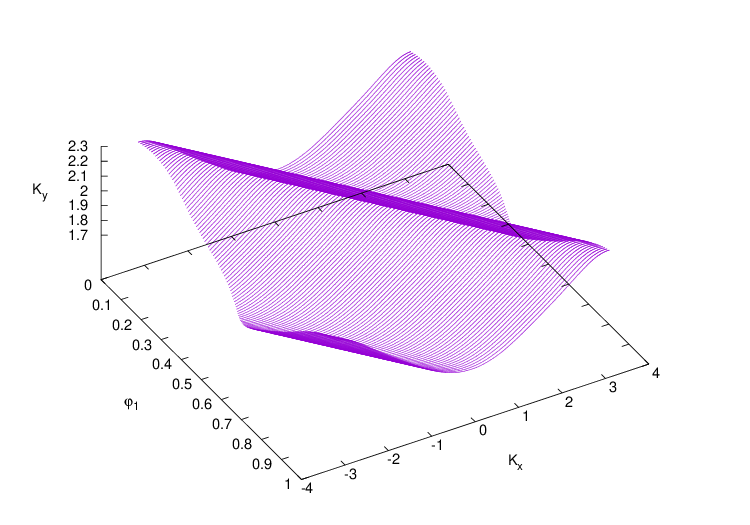}
    \end{minipage}
    \hfill
    \begin{minipage}{0.45\textwidth}
        \centering
        \includegraphics[width=\textwidth]{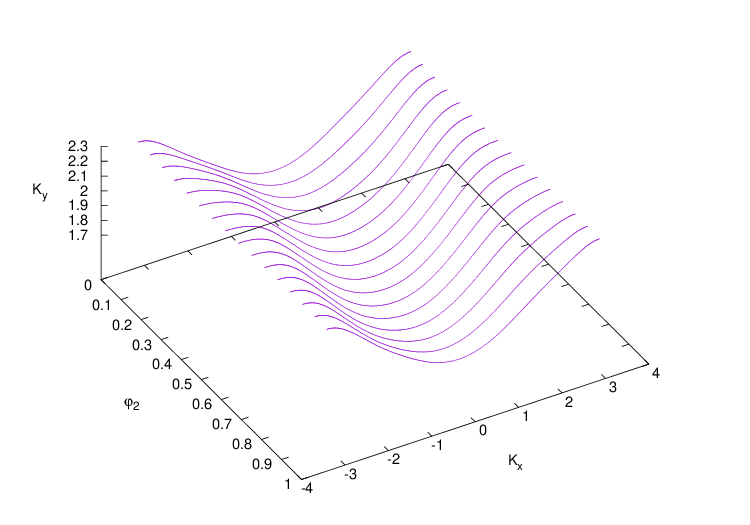}
    \end{minipage}
    \caption{The parameters are $\epsilon_2 = 0.320$ and $\epsilon_3 = 0.155$. The coefficients in each direction are $2^9$, $2^9$, and $2^7$ for $\theta$, $\varphi_1$, and $\varphi_2$, respectively.}
    \label{fig:Torusf1-320}
\end{figure}

In his work, del-Castillo (2000), 
\cite{del2000self}, 
presents \textit{the single-wave model} as a further simplification of the vorticity mixing problem. 
Specifically, 
he introduces a pendulum-like Hamiltonian with a time-dependent amplitude, 
described by,
\begin{equation}\label{eq:Trans:del}
    \mathcal{H}=\sum_{j=1}^{N}\left(\dfrac{p_j^2}{2\Gamma_j} -2\Gamma_j \sqrt{J}\cos(q_j -\theta) \right)-UJ \,,
\end{equation}
with $\Gamma_j$ and $U$ constants, 
and $J$, $\theta$ functions of time.
Following \cite{del2000self, martinez2015self}, we adopt 
\[
2\sqrt{J} = \epsilon_2 + \epsilon_3 \cos(\alpha_1 t), \quad \theta = \alpha_2 t,
\]
Let us consider the model by setting $j=1$, $\Gamma_1=1$ and $U=0$, 
then we obtain
\begin{equation}\label{eq:del}
    \H=\dfrac{p^2}{2} - (\epsilon_2 + \epsilon_3 \cos(\alpha_1 t))\cos(q_j -\alpha_2t)  \,.
\end{equation}
We observe that \eqref{eq:del} coincides with equation \eqref{eq:Transport:mod},
when $\epsilon_1=0$.
Therefore the vector field is
\begin{equation}
    \tilde Z_\H(z,\varphi_1,\varphi_2)= \begin{pmatrix}
        Z_\H(z,\varphi_1,\varphi_2) \\
        \dot \varphi_1 \\
        \dot \varphi_2 \\
    \end{pmatrix}=
    \begin{pmatrix}
        p \\
        (\epsilon_2 + \epsilon_3\cos(\varphi_1))\cos(\varphi_2 -q) \\
        \alpha_1 \\
        \alpha_2 \\
    \end{pmatrix}
\end{equation}
In this case the torus is rotational, 
and the parameterization is of the form
\[
        \tilde K (\theta,( \varphi_1,\varphi_2)) = \begin{pmatrix} K(\theta,(\varphi_1,\varphi_2)) \\ \varphi_1 \\ \varphi_2 \end{pmatrix} =
            \begin{pmatrix} \theta + K_{x}^p(\theta,(\varphi_1,\varphi_2))  \\ K_{y}^p(\theta,(\varphi_1,\varphi_2)) \\ \varphi_1 \\ \varphi_2 \end{pmatrix}  \,.
\]
To construct the initial torus,
we set parameters as $\epsilon_1=\epsilon_2=\epsilon_3=0$,
with the initial $q(0)=0$ and $p(0)=2$. 
This implies that the internal frecuency is $\omega=2$,
while the torus componentes are $K_x(\theta,\varphi_1,\varphi_2)=\theta$ and $K_y(\theta,\varphi_1,\varphi_2)=p_0$.
Under these conditions,
the parameterization of the torus is flat.
For the continuation, 
we keep $\epsilon_1=0$ and $\epsilon_2 = \epsilon_3$, 
up to $\epsilon_2=0.320$ and $\epsilon_3=0.0155$,
figure~\ref{fig:Torusf1-320}.
\subsection{Quasi-periodic pendulum}\label{ssec:Pendulum}
In this section, 
we consier $\epsilon_1=1$ in \eqref{eq:Transport:mod}, i.e., 
\begin{equation}\label{eq:Transport:del}
    \H=\dfrac{p^2}{2} +\cos(q)- (\epsilon_2 + \epsilon_3 \cos(\alpha_1 t))\cos(q_j -\alpha_2t)  \,.
\end{equation}
We explore numerically two different families of tori corresponding to rotational and librational movements, 
i.e., 
families that are topologically distinct.
In both cases, 
we integrate the pendulum vector field with $\epsilon_2=\epsilon_3=0$ to obtain $K_0(\theta) := K(\theta, 0)$. 
Afterward, 
we apply the  method described in subsection~\ref{ssec:inv:torus} to construct the initial torus.
\subsubsection{Rotational motions}\label{sssec:Rotational}
To calculate the internal frequency in this case, 
we follow a procedure similar to that of the Tokamak in Section~\ref{ssec:inv:torus}. 
We do the numerical integration of the vector field,
with initial conditions $q(0)=-\pi$ and $p(0)=2$,
and apply Birkhoff average developed in \cite{das2017quantitative}, resulting in $\omega =2.39509749976984749999$.
The following selection of images shows, 
on the left, 
the components of the torus $K$ with respect to the coordinate $\varphi_1$, 
and on the right, 
the same components with respect to $\varphi_2$.
\begin{figure}[htbp]
    \begin{minipage}{0.45\textwidth}
        \centering
        \includegraphics[width=\textwidth]{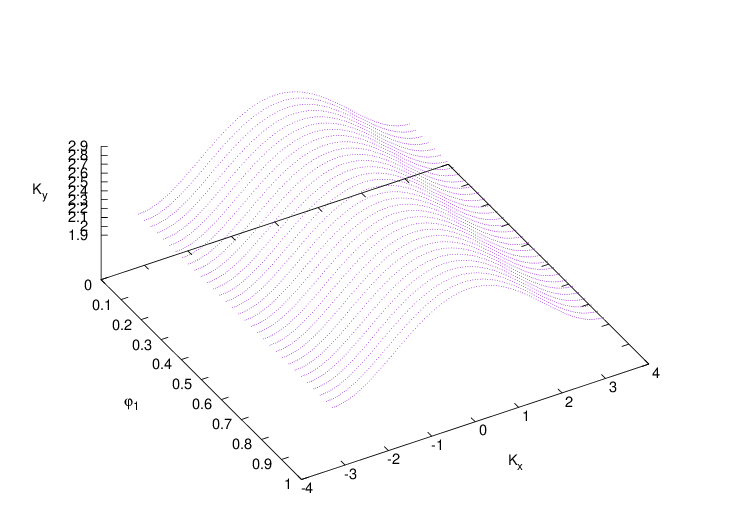}
    \end{minipage}
    \hfill
    \begin{minipage}{0.45\textwidth}
        \centering
        \includegraphics[width=\textwidth]{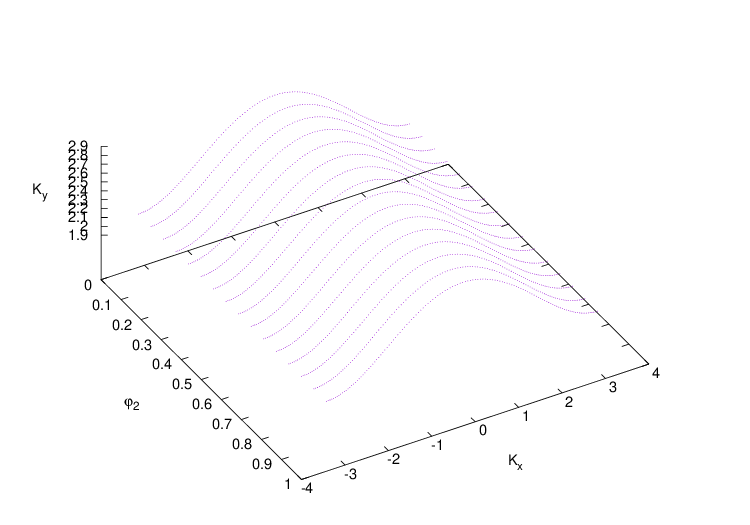}
    \end{minipage}
    \caption{The parameters are $\epsilon_1=1.000$ $\epsilon_2 = 0.000$ and $\epsilon_3 = 0.000$. The coefficients in each direction are $2^7$, $2^5$, and $2^4$ for $\theta$, $\varphi_1$, and $\varphi_2$, respectively.}
\end{figure}

\begin{figure}[htbp]
    \begin{minipage}{0.45\textwidth}
        \centering
        \includegraphics[width=\textwidth]{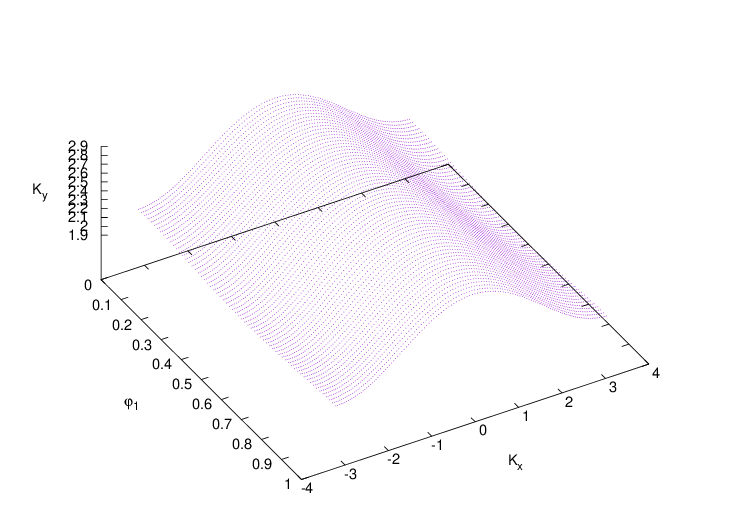}
    \end{minipage}
    \hfill
    \begin{minipage}{0.45\textwidth}
        \centering
        \includegraphics[width=\textwidth]{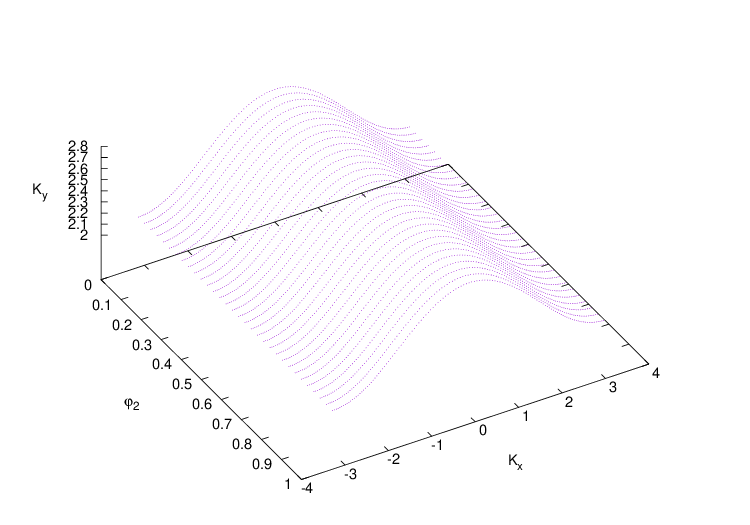}
    \end{minipage}
    \caption{The parameters are $\epsilon_1=0.096$ $\epsilon_2 = 0.040$ and $\epsilon_3 = 0.040$. The coefficients in each direction are $2^7$, $2^6$, and $2^5$ for $\theta$, $\varphi_1$, and $\varphi_2$, respectively.}
\end{figure}

\begin{figure}[htbp]
    \begin{minipage}{0.45\textwidth}
        \centering
        \includegraphics[width=\textwidth]{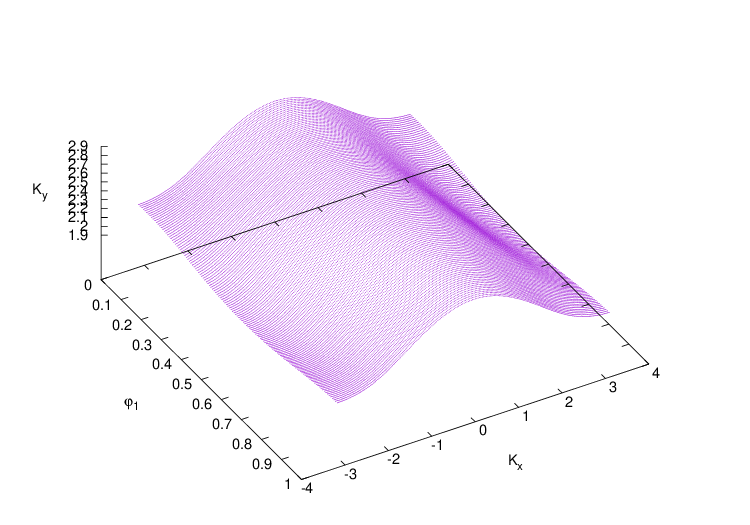}
    \end{minipage}
    \hfill
    \begin{minipage}{0.45\textwidth}
        \centering
        \includegraphics[width=\textwidth]{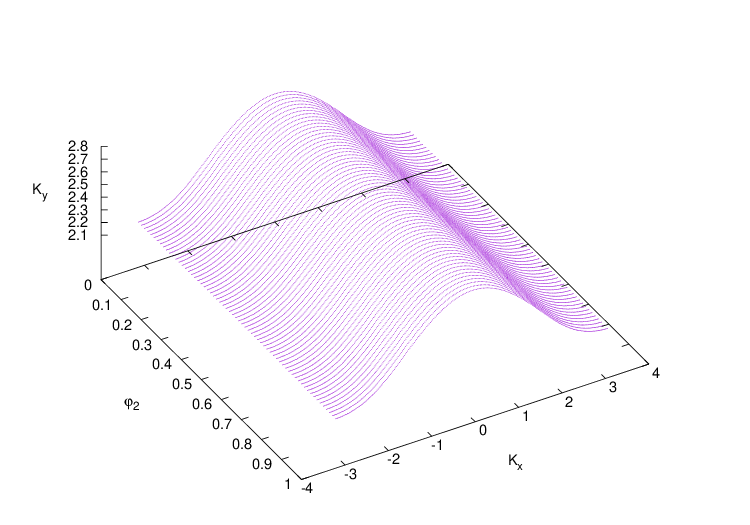}
    \end{minipage}
    \caption{The parameters are $\epsilon_1=0.092$ $\epsilon_2 = 0.080$ and $\epsilon_3 = 0.040$. The coefficients in each direction are $2^8$, $2^5$, and $2^4$ for $\theta$, $\varphi_1$, and $\varphi_2$, respectively.}
\end{figure}
\break
\subsubsection{Librational motions}\label{sssec:Librational}
For the librational tori, 
we consider the initial conditions $q(0)=\frac{\pi}{2}$ and $p(0)=0$. 
From these, 
using Jacobi elliptic integrals, 
internal frequency is $\omega =0.84721308479397908659$.
Similar to the previous section, 
the following set of images provides a comparative view of the torus $K$. 
On the left, the components are plotted with respect to the coordinate $\varphi_1$, 
while on the right, the same components are shown relative to $\varphi_2$.
\begin{figure}[htbp]
    \centering
    \begin{minipage}{0.45\textwidth}
        \centering
        \includegraphics[width=\textwidth]{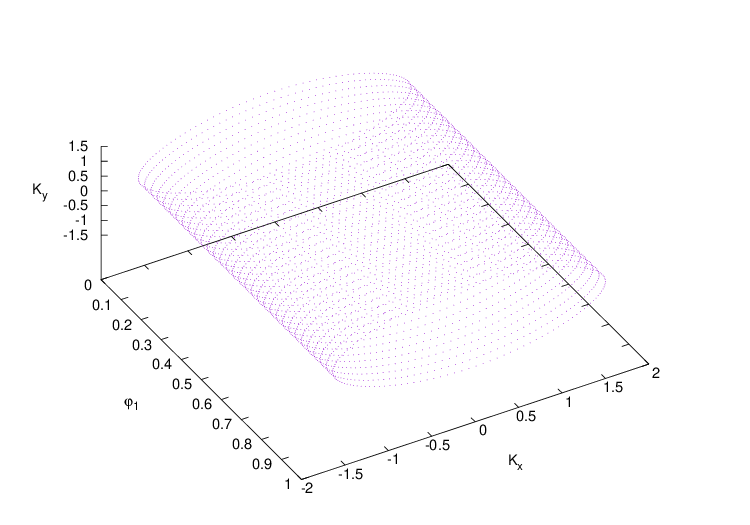}
    \end{minipage}
    \hfill
    \begin{minipage}{0.45\textwidth}
        \centering
        \includegraphics[width=\textwidth]{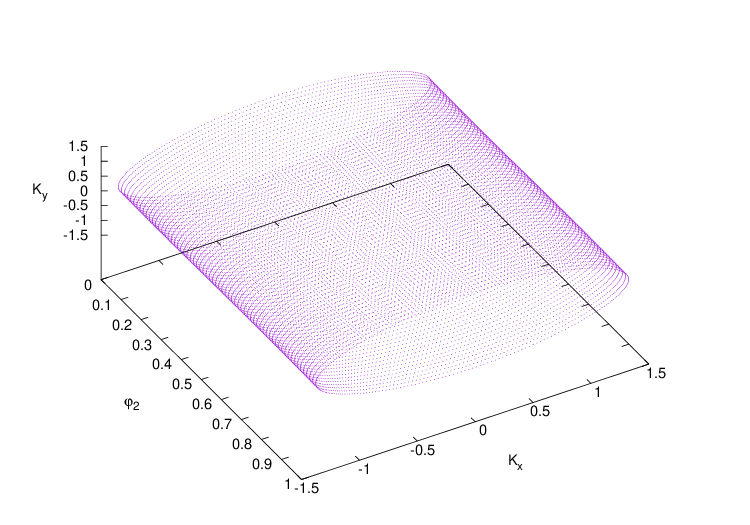}
    \end{minipage}
    \caption{The parameters are $\epsilon_1 = 1.000$, $\epsilon_2 = 0.000$ and $\epsilon_3 = 0.000$. The coefficients in each direction are $2^7$, $2^5$, and $2^4$ for $\theta$, $\varphi_1$, and $\varphi_2$, respectively.}
    \vspace{5mm}
\end{figure}

\begin{figure}[htbp]
    \begin{minipage}{0.45\textwidth}
        \centering
        \includegraphics[width=\textwidth]{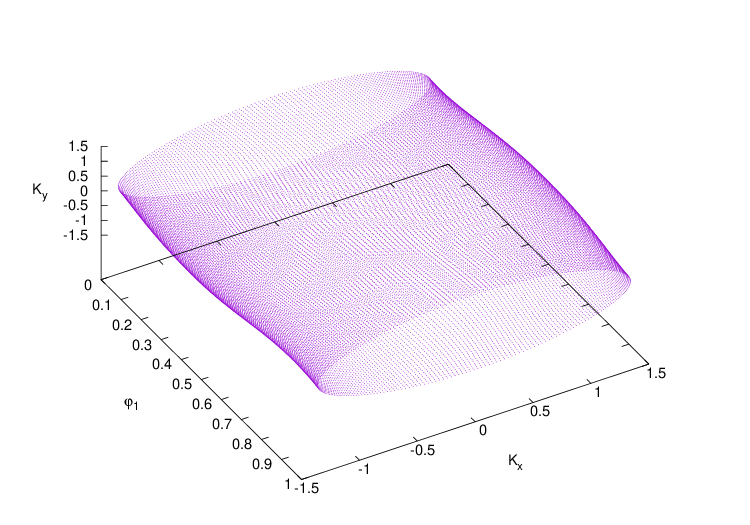}
    \end{minipage}
    \hfill
    \begin{minipage}{0.45\textwidth}
        \centering
        \includegraphics[width=\textwidth]{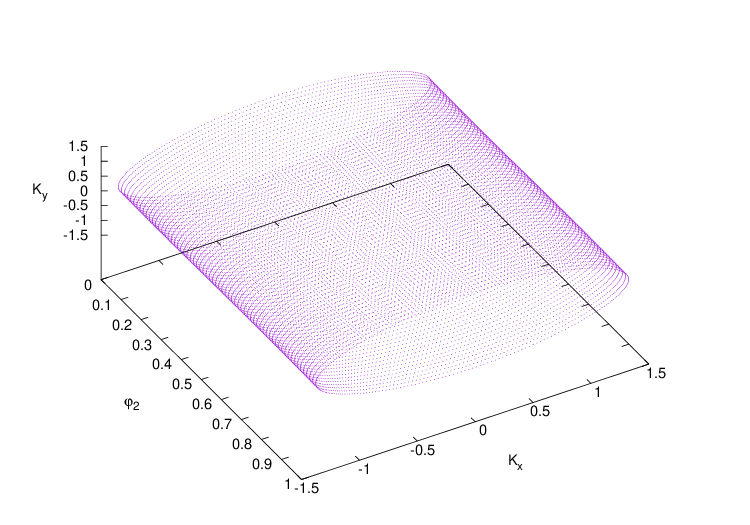}
    \end{minipage}
    \caption{The parameters are $\epsilon_1=0.930$ $\epsilon_2 = 0.035$ and $\epsilon_3 = 0.007$. The coefficients in each direction are $2^8$, $2^7$, and $2^6$ for $\theta$, $\varphi_1$, and $\varphi_2$, respectively.}
    \vspace{5mm}
\end{figure}

\begin{figure}[htbp]
    \begin{minipage}{0.45\textwidth}
        \centering
        \includegraphics[width=\textwidth]{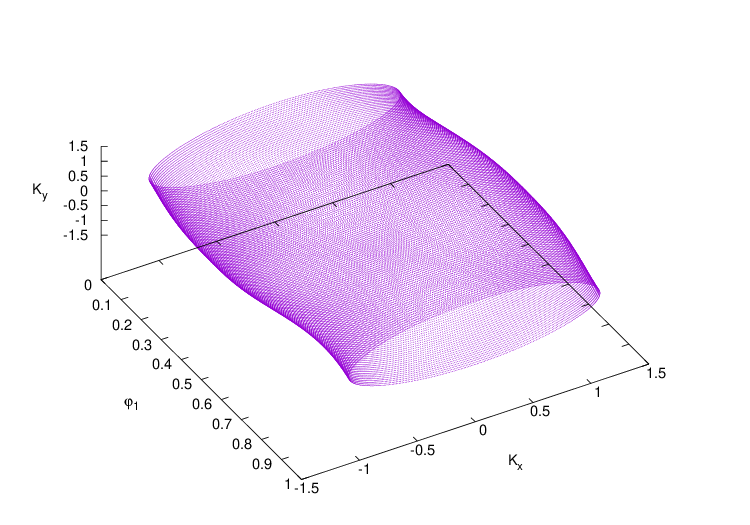}
    \end{minipage}
    \hfill
    \begin{minipage}{0.45\textwidth}
        \centering
        \includegraphics[width=\textwidth]{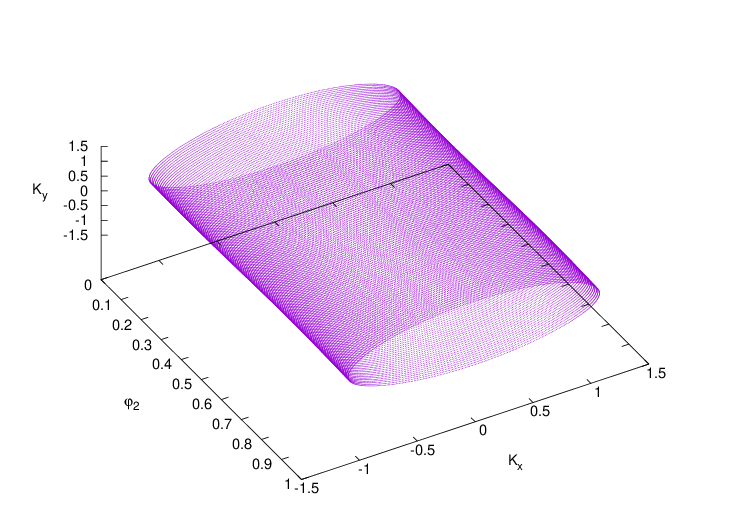}
    \end{minipage}
    \caption{The parameters are $\epsilon_1=0.860$ $\epsilon_2 = 0.070$ and $\epsilon_3 = 0.014$. The coefficients in each direction are $2^9$, $2^8$, and $2^7$ for $\theta$, $\varphi_1$, and $\varphi_2$, respectively.}
\end{figure}
\break 
\subsection{Computation details}
The numerical simulations were executed on a computer equipped with an \textbf{Intel(R) Xeon(R) E5-2665 CPU}, 
operating at \textbf{2.40 GHz} with a maximum boost frequency of \textbf{3.10 GHz}. 
This processor features a \textbf{16-core} architecture, 
allowing for \textbf{16 threads} due to its hyper-threading capability, 
with \textbf{8 cores per socket}. 
The system architecture is \textbf{x86\_64}, 
supporting both \textbf{32-bit and 64-bit} modes. 
The CPU's cache hierarchy includes \textbf{256 KB} of L1 data cache per core, 
\textbf{256 KB} of L1 instruction cache per core, 
\textbf{2 MB} of L2 cache per core, 
and \textbf{20 MB} of shared L3 cache. 
The machine operates on a single NUMA node with all CPUs (0-15) accessible, 
and it supports virtualization through \textbf{VT-x} technology. 
The system is running on a \textbf{Linux} operating system. 
The machine is equipped with \textbf{32 GB} of total RAM, 
Additionally, 
the system has a \textbf{32 GB} swap space, 
The calculations were performed using \textbf{C++} with \texttt{long double},
providing a precision of up to 20 significant digits.

\section{Conclusions}

In this work, 
we have presented two applications of the  method for the existence and stability of invariant tori in time-dependent quasi-periodic Hamiltonian systems. 
Unlike other approaches, 
our method is applied directly to the vector field of the system, 
avoiding the need to construct a quasi-periodic map, 
as is commonly done in studies based on the standard map. 
This formulation allows for a more direct and flexible implementation in different physical contexts, 
facilitating its adaptation to problems with various dynamical structures.

Additionally, 
we exploit the fibered structure of the system to reduce the problem’s dimension from $2n + \ell$ to $2n$, 
where the external frequencies are naturally incorporated into the formulation. 
This not only simplifies the theoretical analysis but also improves computational efficiency. 
In particular, 
the torsion expression used in this work avoids the computation of the Lie derivative, 
as required in other methods, 
significantly reducing computational complexity. 
Instead, 
our method only requires evaluating the torus parameterization, 
in the derivative of the vector field, 
and the product of geometric objects, 
leading to faster and more precise calculations. 
This reduction in computational cost is especially relevant for applications that require high numerical precision.

Furthermore, 
our approach enables more efficient Computer-Assisted Proofs (CAPs) compared to methods based on the system’s flow. 
Since we work directly with the vector field rather than integrating the dynamics, 
we avoid the need to validate numerical integration through a CAP, 
which can be computationally expensive and numerically delicate. 
This aspect makes our method particularly suitable for applications requiring rigorous validation with lower computational effort.

In cases where the method converges, 
the number of Fourier nodes doubles at each iteration, 
leading to an increasing demand for computational resources. 
In the case of the Tokamak, the expansion is first attempted in the $\theta$ direction; 
if this does not work, 
it is tried in $\varphi$, 
and if it still does not converge, 
it is performed in both directions simultaneously. 
A similar procedure is followed for coupled pendulums. 
This behavior makes RAM management a crucial aspect of the method’s implementation. 
In the future, when Computer-Assisted Proofs are carried out, 
it will be essential to perform calculations in multiple precision, 
further increasing memory demand and requiring advanced resource management strategies. 
This challenge represents an interesting problem that is worth addressing in future research.

Our results strengthen the applicability of KAM methods in time-dependent Hamiltonian systems and open new perspectives for their implementation in concrete physical problems, 
such as the dynamics of confined plasmas and geometric transport models. 
Future work could explore extending this approach to systems with a larger number of degrees of freedom and integrating it with advanced numerical techniques for high-dimensional problems. Additionally, 
combining this method with computational analysis approaches and computer-assisted proofs could provide new tools for the rigorous study of the stability and persistence of invariant tori in more general scenarios.

The method is designed to compute tori of any dimension.
The only restrintion comes from the computational resources.

\section{Acknowledgements}
We would like to express our gratitude to the following organizations for their support: DGAPA-UNAM through project PAPIIT IN103423, IN104725
CONACYT for the PhD fellowship, 
and the support received from the project PID2021-125535NB-I00 (MCIU/AEI/FEDER, UE). 
Additionally, A. H. acknowledge the funding received from the Severo Ochoa and María de Maeztu Program for Centers and Units of Excellence in R\&D (CEX2020-001084-M).
We are also grateful to
D. del Castillo-Negrete, D. Mart\'inez-del R\'io, A. Olvera, A. Vieiro,
J.-Ll. Figueras, J. M.  Mondelo and \'A. Fern\'andez-Mora for fruitful discussions.
R.C. 
and P.P. would like to express their sincere gratitude to the Departments of Mathematics and Computer Science at the University of Barcelona.
Additionally, 
P.P. would also like to thank the Department of Mathematics at Uppsala University, 
for their warm hospitality during my stay at their institutions while carrying out this work.\\

\bibliographystyle{alpha}

\begin{thebibliography}{MdRdCNOC15}

\bibitem[Abd06]{abdullaev2006construction}
Sadrilla~S Abdullaev.
\newblock {\em Construction of mappings for Hamiltonian systems and their
  applications}, volume 691.
\newblock Springer, 2006.

\bibitem[Arn63]{Arnold63a}
V.I. Arnold.
\newblock Proof of a theorem of {A}. {N}. {K}olmogorov on the preservation of
  conditionally periodic motions under a small perturbation of the
  {H}amiltonian.
\newblock {\em Uspehi Mat. Nauk}, 18(5 (113)):13--40, 1963.

\bibitem[Art72]{artsimovich1972tokamak}
LA~Artsimovich.
\newblock Tokamak devices.
\newblock {\em Nuclear Fusion}, 12(2):215, 1972.

\bibitem[BMO09]{barrabes2009numerical}
Esther Barrabes, Josep~M Mondelo, and Merce Olle.
\newblock Numerical continuation of families of homoclinic connections of
  periodic orbits in the rtbp.
\newblock {\em Nonlinearity}, 22(12):2901, 2009.

\bibitem[Boo04]{boozer2004physics}
Allen~H Boozer.
\newblock Physics of magnetically confined plasmas.
\newblock {\em Reviews of modern physics}, 76(4):1071--1141, 2004.

\bibitem[CC06]{celletti2006kam}
Alessandra Celletti and Luigi Chierchia.
\newblock {KAM} tori for n-body problems: a brief history.
\newblock In {\em Periodic, Quasi-Periodic and Chaotic Motions in Celestial
  Mechanics: Theory and Applications: Selected papers from the Fourth Meeting
  on Celestial Mechanics, CELMEC IV San Martino al Cimino (Italy), 11--16
  September 2005}, pages 117--139. Springer, 2006.

\bibitem[CCGdlL24]{calleja2024accurate}
Renato Calleja, Alessandra Celletti, Joan Gimeno, and Rafael de~la Llave.
\newblock Accurate computations up to breakdown of quasi-periodic attractors in
  the dissipative spin--orbit problem.
\newblock {\em Journal of Nonlinear Science}, 34(1):12, 2024.

\bibitem[CDlL10]{calleja2010numerically}
Renato Calleja and Rafael De~la Llave.
\newblock A numerically accessible criterion for the breakdown of
  quasi-periodic solutions and its rigorous justification.
\newblock {\em Nonlinearity}, 23(9):2029, 2010.

\bibitem[CHP25]{calleja2025constructive}
Renato Calleja, Alex Haro, and Pedro Porras.
\newblock Constructive approaches to qp-time-dependent kam theory for
  lagrangian tori in hamiltonian systems.
\newblock {\em arXiv preprint arXiv:2503.09740}, 2025.

\bibitem[CL10]{calleja2010computation}
Renato Calleja and Rafael Llave.
\newblock Computation of the breakdown of analyticity in statistical mechanics
  models: numerical results and a renormalization group explanation.
\newblock {\em Journal of Statistical Physics}, 141(6):940, 2010.

\bibitem[CVC{\etalchar{+}}05]{chandre2005control}
Cristel Chandre, Michel Vittot, Guido Ciraolo, Ph~Ghendrih, and Ricardo Lima.
\newblock Control of stochasticity in magnetic field lines.
\newblock {\em Nuclear Fusion}, 46(1):33, 2005.

\bibitem[dCN98a]{del1998weakly}
D~del Castillo-Negrete.
\newblock Weakly nonlinear dynamics of electrostatic perturbations in
  marginally stable plasmas.
\newblock {\em Physics of Plasmas}, 5(11):3886--3900, 1998.

\bibitem[dCN98b]{del1998nonlinear}
Diego del Castillo-Negrete.
\newblock Nonlinear evolution of perturbations in marginally stable plasmas.
\newblock {\em Physics Letters A}, 241(1-2):99--104, 1998.

\bibitem[dCN00]{del2000self}
Diego del Castillo-Negrete.
\newblock Self-consistent chaotic transport in fluids and plasmas.
\newblock {\em Chaos: An Interdisciplinary Journal of Nonlinear Science},
  10(1):75--88, 2000.

\bibitem[DlL{\etalchar{+}}01]{de2001tutorial}
Rafael De~la Llave et~al.
\newblock A tutorial on {KAM} theory.
\newblock In {\em Proceedings of Symposia in Pure Mathematics}, volume~69,
  pages 175--296. Providence, RI; American Mathematical Society; 1998, 2001.

\bibitem[dlLGJV05]{de2005kam}
Rafael de~la Llave, Alejandra Gonz{\'a}lez, {\`A}ngel Jorba, and Jordi
  Villanueva.
\newblock {KAM} theory without action-angle variables.
\newblock {\em Nonlinearity}, 18(2):855, 2005.

\bibitem[DSSY17]{das2017quantitative}
Suddhasattwa Das, Yoshitaka Saiki, Evelyn Sander, and James~A Yorke.
\newblock Quantitative quasiperiodicity.
\newblock {\em Nonlinearity}, 30(11):4111, 2017.

\bibitem[Dum14]{dumas2014kam}
H~Scott Dumas.
\newblock {\em Kam Story, The: A Friendly Introduction To The Content, History,
  And Significance Of Classical Kolmogorov-arnold-moser Theory}.
\newblock World Scientific Publishing Company, 2014.

\bibitem[FH25]{figueras2025sun}
Jordi-Llu{\'\i}s Figueras and Alex Haro.
\newblock Sun--jupiter--saturn system may exist: A verified computation of
  quasiperiodic solutions for the planar three-body problem.
\newblock {\em Journal of Nonlinear Science}, 35(1):1--20, 2025.

\bibitem[FHL17]{FiguerasHL17}
J.-Ll. Figueras, A.~Haro, and A.~Luque.
\newblock Rigorous {C}omputer-{A}ssisted {A}pplication of {KAM} {T}heory: {A}
  {M}odern {A}pproach.
\newblock {\em Found. Comput. Math.}, 17(5):1123--1193, 2017.

\bibitem[FMHM24]{fernandez2024flow}
{\'A}lvaro Fern{\'a}ndez-Mora, Alex Haro, and Josep-Maria Mondelo.
\newblock Flow map parameterization methods for invariant tori in
  quasi-periodic hamiltonian systems.
\newblock {\em SIAM Journal on Applied Dynamical Systems}, 23(1):127--166,
  2024.

\bibitem[GHdlL22]{gonzalez2022efficient}
Alejandra Gonz{\'a}lez, {\`A}lex Haro, and Rafael de~la Llave.
\newblock Efficient and reliable algorithms for the computation of non-twist
  invariant circles.
\newblock {\em Foundations of Computational Mathematics}, 22(3):791--847, 2022.

\bibitem[GM01]{gomez2001dynamics}
Gerard G{\'o}mez and JOSEP~M Mondelo.
\newblock The dynamics around the collinear equilibrium points of the rtbp.
\newblock {\em Physica D: Nonlinear Phenomena}, 157(4):283--321, 2001.

\bibitem[HCF{\etalchar{+}}16]{HaroCFLM16}
\`A. Haro, M.~Canadell, J.-Ll. Figueras, A.~Luque, and J.-M. Mondelo.
\newblock {\em The parameterization method for invariant manifolds}, volume 195
  of {\em Applied Mathematical Sciences}.
\newblock Springer, [Cham], 2016.
\newblock From rigorous results to effective computations.

\bibitem[HL19]{haro2019posteriori}
Alex Haro and Alejandro Luque.
\newblock A-posteriori {KAM} theory with optimal estimates for partially
  integrable systems.
\newblock {\em Journal of Differential Equations}, 266(2-3):1605--1674, 2019.

\bibitem[Kol54]{kolmogorov1954conservation}
Andrey~Nikolaevich Kolmogorov.
\newblock On conservation of conditionally periodic motions for a small change
  in hamilton's function.
\newblock In {\em Dokl. Akad. Nauk SSSR}, volume~98, pages 527--530, 1954.

\bibitem[MdRdCNOC15]{martinez2015self}
D~Mart{\'\i}nez-del R{\'\i}o, Diego del Castillo-Negrete, Arturo Olvera, and
  Renato Calleja.
\newblock Self-consistent chaotic transport in a high-dimensional mean-field
  hamiltonian map model.
\newblock {\em Qualitative theory of dynamical systems}, 14(2):313--335, 2015.

\bibitem[MMP84]{mackay1984transport}
RS~MacKay, JD~Meiss, and IC~Percival.
\newblock Transport in hamiltonian systems.
\newblock {\em Physica D: Nonlinear Phenomena}, 13(1-2):55--81, 1984.

\bibitem[M{\"o}s62]{moser1962invariant}
J~M{\"o}ser.
\newblock On invariant curves of area-preserving mappings of an annulus.
\newblock {\em Nachr. Akad. Wiss. G{\"o}ttingen, II}, pages 1--20, 1962.

\bibitem[VL21]{valvo2021hamiltonian}
Lorenzo Valvo and Ugo Locatelli.
\newblock Hamiltonian control of magnetic field lines: Computer assisted.
\newblock 2021.

\end{thebibliography}

\newcommand{\etalchar}[1]{$^{#1}$}

\break
\appendix
 \section{Table of time and memory for Quasi-Newton method to Tokamak}\label{sec:appen:Tokamak}
The following table shows the time and memory usage percentage at each step of the iterations of the quasi-Newton method. 
The steps are as follows: 
1. Construction of the tangent frame, 
2. Construction of torsion, 
3. Correction in the symplectic frame, 
and 4. Computation of the new parameterization.\\
\begin{tabular}{|c|c|c|c|c|}
    \hline
    \multicolumn{5}{|c|}{\textbf{Iteration 1}} \\
    \hline
    \hline
    Step & Execution Time & \ Memory Usage &  \multicolumn{2}{c}{Objects} \vline  \\ 
    \hline
    1 & \texttt{26.9559} (s) & \texttt{2.0} \% & $\norm{L}=$\texttt{2.602235e+00} & $\norm{\Omega_L}=$\texttt{0.000000}  \\ 
    \hline
    2 & \texttt{43.6666} (s) & \texttt{3.2} \% & $\norm{\Lie{}L}=$\texttt{0.32544} & $\norm{T}=$\texttt{1.15074036e+01} \\ 
    \hline
    3 & \texttt{18.5556} (s) & \texttt{3.2} \% &  $\aver{T}=$\texttt{-2.35551120e+00}&  $\aver{T^{-1}}=$\texttt{-4.24536295e-01}  \\ 
    \hline
    4 & \texttt{0.776569} (s) & \texttt{3.2} \% & $\norm{\Delta K}=$\texttt{8.997e-03} & $\norm{E}=$\texttt{2.71804e-05}  \\ 
    \hline
    \hline
    \multicolumn{5}{|c|}{\textbf{Iteration 2}} \\
    \hline
    \hline
    Step & Execution Time & \ Memory Usage &  \multicolumn{2}{c}{Objects} \vline  \\ 
    \hline
    1 & \texttt{26.9567} (s) & \texttt{2.0} \% & $\norm{L}=$\texttt{2.757424e+00} & $\norm{\Omega_L}=$\texttt{0.000000}  \\ 
    \hline
    2 & \texttt{43.6666} (s) & \texttt{3.2} \% & $\norm{\Lie{}L}=$\texttt{0.322738} & $\norm{T}=$\texttt{1.13239957e+01} \\ 
    \hline
    3 & \texttt{18.5332} (s) & \texttt{3.2} \% &  $\aver{T}=$\texttt{-2.35326725e+00}&  $\aver{T^{-1}}=$\texttt{-4.24941112e-01}  \\ 
    \hline
    4 & \texttt{0.775557} (s) & \texttt{3.2} \% & $\norm{\Delta K}=$\texttt{7.180e-04} & $\norm{E}=$\texttt{5.35071e-06}  \\ 
    \hline
    \hline
    \multicolumn{5}{|c|}{\textbf{Iteration 3}} \\
    \hline
    \hline
    Step & Execution Time & \ Memory Usage &  \multicolumn{2}{c}{Objects} \vline  \\ 
    \hline
    1 & \texttt{26.978} (s) & \texttt{2.0} \% & $\norm{L}=$\texttt{2.789292e+00} & $\norm{\Omega_L}=$\texttt{0.000000}  \\ 
    \hline
    2 & \texttt{43.7547} (s) & \texttt{3.2} \% & $\norm{\Lie{}L}=$\texttt{0.325661} & $\norm{T}=$\texttt{1.13989675e+01} \\ 
    \hline
    3 & \texttt{18.5258} (s) & \texttt{3.2} \% &  $\aver{T}=$\texttt{-2.35312392e+00}&  $\aver{T^{-1}}=$\texttt{-4.24966994e-01}  \\ 
    \hline
    4 & \texttt{0.780214} (s) & \texttt{3.2} \% & $\norm{\Delta K}=$\texttt{1.940e-05} & $\norm{E}=$\texttt{4.22474e-08}  \\ 
    \hline
    \hline
    \multicolumn{5}{|c|}{\textbf{Iteration 4}} \\
    \hline
    \hline
    Step & Execution Time & \ Memory Usage &  \multicolumn{2}{c}{Objects} \vline  \\ 
    \hline
    1 & \texttt{27.0355} (s) & \texttt{2.0} \% & $\norm{L}=$\texttt{2.785728e+00} & $\norm{\Omega_L}=$\texttt{0.000000}  \\ 
    \hline
    2 & \texttt{43.8932} (s) & \texttt{3.2} \% & $\norm{\Lie{}L}=$\texttt{0.325439} & $\norm{T}=$\texttt{1.13621447e+01} \\ 
    \hline
    3 & \texttt{18.5478} (s) & \texttt{3.2} \% &  $\aver{T}=$\texttt{-2.35313653e+00}&  $\aver{T^{-1}}=$\texttt{-4.24964718e-01}  \\ 
    \hline
    4 & \texttt{0.778876} (s) & \texttt{3.2} \% & $\norm{\Delta K}=$\texttt{1.895e-08} & $\norm{E}=$\texttt{2.70577e-11}  \\ 
    \hline
    \hline
    \multicolumn{5}{|c|}{\textbf{Iteration 5}} \\
    \hline
    \hline
    Step & Execution Time & \ Memory Usage &  \multicolumn{2}{c}{Objects} \vline  \\ 
    \hline
    1 & \texttt{26.9686} (s) & \texttt{2.0} \% & $\norm{L}=$\texttt{2.785741e+00} & $\norm{\Omega_L}=$\texttt{0.000000}  \\ 
    \hline
    2 & \texttt{56.311} (s) & \texttt{3.2} \% & $\norm{\Lie{}L}=$\texttt{0.32544} & $\norm{T}=$\texttt{1.136222460e+01} \\ 
    \hline
    3 & \texttt{18.4804} (s) & \texttt{3.2} \% &  $\aver{T}=$\texttt{-2.35313653e+00}&  $\aver{T^{-1}}=$\texttt{-4.24964718e-01}  \\ 
    \hline
    4 & \texttt{0.778876} (s) & \texttt{3.2} \% & $\norm{\Delta K}=$\texttt{2.706e-08} & $\norm{E}=$\texttt{9.63249e-16}  \\ 
    \hline

\end{tabular}

\end{document}